\documentclass[a4paper,fleqn]{cas-sc}
\usepackage{subcaption}
\usepackage[numbers,sort&compress]{natbib}
\newcommand{\vect}[1]{\boldsymbol{\mathbf{#1}}}
\usepackage{cleveref}
 \usepackage{amsthm}
  \newtheorem{theorem}{Theorem}
 
\newdefinition{remark}{Remark}

\def\tsc#1{\csdef{#1}{\textsc{\lowercase{#1}}\xspace}}
\tsc{WGM}
\tsc{QE}
\tsc{EP}
\tsc{PMS}
\tsc{BEC}
\tsc{DE}

\begin{document}
\let\WriteBookmarks\relax
\def\floatpagepagefraction{1}
\def\textpagefraction{.001}
\shorttitle{Weight-Adjusted Entropy-Stable BR2 and CDG discretizations with Entropy-Stable BCs}
\shortauthors{Pranshul Thakur and Siva Nadarajah}

\title [mode = title]{Weight-Adjusted Entropy-Stable BR2 and CDG Viscous Discretizations with Wall Boundary Conditions for Compressible Flows on Curvilinear Grids}                      
\tnotemark[1]

\tnotetext[1]{This document is the result of the research
   project funded by the Fonds de recherche du Québec – Nature et technologies (FRQNT) Doctoral Award, McGill Engineering Doctoral Award (MEDA) and the Natural Sciences and Engineering Research Council of Canada (NSERC) Discovery Grant RGPIN-2019-04791.}

\author[1]{Pranshul Thakur}
\cormark[1]
\ead{pranshul.thakur@mail.mcgill.ca}

\credit{Methodology, Conceptualization, Software, Writing – original draft}

\affiliation[1]{organization={Department of Mechanical Engineering, McGill University}, 
                city={Montréal},
                postcode={H3A 0G4}, 
                state={Québec},
                country={Canada}}

\author[1]{Siva Nadarajah}
\ead{siva.nadarajah@mcgill.ca}

\credit{Methodology, Funding acquisition, Software, Writing - Original draft}

\cortext[cor2]{Corresponding author}

\begin{abstract}
Entropy-stable discontinuous Galerkin schemes are known for their robustness. These schemes satisfy a semi-discrete entropy inequality, mimicking the entropy inequality of the continuous governing equations. Entropy-stable viscous discretizations with entropy-stable wall boundary conditions have been used in the literature, using a scheme similar to BR1. On the other hand, viscous discretizations with compact stencil, such as BR2 and CDG schemes, are more efficient for computation and for forming the Jacobian. Moreover, it is known that the weight-adjusted formulations have better efficiency with regard to the computational and storage costs. To that end, we derive the weight-adjusted entropy-stable BR2 and CDG schemes with modifications for the numerical flux at the boundary to account for the entropy-stable boundary conditions.
\end{abstract}




\maketitle

\section{Introduction}
High-order numerical schemes have been gaining traction due to their improved accuracy, at computational costs similar to the first- and second-order schemes. Discontinuous Galerkin (DG) scheme is one of the most popular high-order schemes, which is known for its compact stencil and better parallelizability \cite{hesthaven2008nodal}. Compact stencil is an important property, as it allows better parallelizability and ease of computing the Jacobian. The Jacobian of a numerical scheme is useful for multiple purposes, such as to run an implicit time-stepping scheme, for aerodynamic shape optimization and for adjoint-based mesh adaptation. 

Since high-order methods add less dissipation, they have issues related to stability and robustness in the presence of shocks and under-resolved features. Several methods have been used in the past to improve the robustness of high-order schemes, such as artificial viscosity \cite{subcell_shock_capturing}, filtering \cite{filtering}, and entropy stability \cite{alex_nsfr,chan_entropystable_dg}. Out of these approaches, the entropy-stable schemes seem promising due to their ability to achieve robustness without adding additional dissipation. It is well-known that several governing equations have a scalar quantity, known as entropy, which satisfies a scalar inequality \cite{lax_entropy}. While conventional discretization schemes do not necessarily satisfy the entropy inequality discretely, entropy-stable schemes ensure that the scalar inequality is discretely satisfied. For linear advection and Burgers equations, the entropy function is a Lyapunov function, and entropy stability ensures nonlinear stability of the scheme. For compressible Euler and Navier-Stokes equations, while there does not exist a formal proof for the boundedness of the solution, several studies \cite{chan_entropystable_dg,alex_nsfr,julien_les} have shown that the entropy-stable discretizations lead to a robust scheme, for which the solution remains bounded even for the test cases involving shocks and under-resolved features. 

Entropy-stable schemes have been extended to discretize the viscous numerical fluxes \cite{May_viscous_entropystableDG,viscous_sipg_br1}. In \cite{May_viscous_entropystableDG}, the authors discussed the entropy-stable treatment for the volume and interface terms for the Symmetric Interior Penalty (SIP), BR2 and Local discontinuous Galerkin (LDG) schemes, yielding entropy stability for periodic boundary conditions. In \cite{viscous_sipg_br1}, the authors developed a SIPG-BR1 hybrid scheme, as well as a BR2 formulation that was entropy stable for periodic boundary conditions. However, for non-periodic boundary conditions, we note that additional modifications are required for the entropy-stable imposition of the boundary conditions. Recent studies in the literature have addressed entropy-stable imposition of boundary conditions for compressible Euler and Navier-Stokes equations \cite{chan_entropystable_viscous,svard_entropystablebc1,svard_entropystablebc2}. In \cite{chan_entropystable_viscous}, the authors introduced a modified entropy-stable BR1 discontinuous Galerkin scheme and provided the entropy-stable boundary conditions for several viscous wall boundary conditions: symmetric, no-slip, isothermal and adiabatic. They proved that the continuous entropy inequality, including the boundary conditions, is discretely satisfied with these modifications. 

However, we note that schemes such as BR1 and LDG do not have a compact stencil. Viscous discretization schemes with a compact stencil are useful because, in such cases, it is easier to form the Jacobian. The Jacobian of the discretization is required for implicit time advancement as well as for optimization and mesh adaptation. In the current work, we extend the entropy-stable boundary treatment in \cite{chan_entropystable_viscous} to two viscous discretization schemes with a compact stencil: BR2 and CDG. We note that the entropy-stable discretizations of BR2 have been investigated in the past in \cite{May_viscous_entropystableDG} and \cite{viscous_sipg_br1}. However, the approach in \cite{viscous_sipg_br1} is only entropy-stable for orthogonal grids, while both the approaches in \cite{May_viscous_entropystableDG} and \cite{viscous_sipg_br1} use standard numerical fluxes at the boundary, which are not entropy-stable for the boundary terms. In addition to the entropy-stable boundary conditions, another novelty of the current approach is that it uses a weight-adjusted formulation for curvilinear grids. The weight-adjusted formulation is considered to be computationally efficient as it can be used with sum-factorization and has low storage requirements \cite{chan_weightadjusted_jcp}.     

We proceed in the following order: in Section 2, we review the entropy stability of the compressible Navier-Stokes equations, along with the continuous entropy inequality and the boundary condition satisfied by the entropy variable. In Section 3, we provide the notations and give an overview of the weight-adjusted formulation. In Section 4, we derive the weight-adjusted entropy-stable formulation for a general viscous discretization, which is applicable for curvilinear grids. In Section 4.1, we provide the specific numerical fluxes for the weight-adjusted BR2 scheme, as well as the modifications to the boundary numerical flux for entropy stability. We refer to this discretization as BR2-E. In Section 4.2, we provide boundary and interior-face numerical fluxes and prove entropy stability for the weight-adjusted CDG-E viscous discretization. Finally, in Section 5, we present two- and three-dimensional test cases to test the entropy stability and convergence orders of the schemes.

\section{Entropy stability of the Navier-Stokes equations}
\label{sec:entropy_stable_ns}

The Navier-Stokes equations in conservative form can be written as:
\begin{equation}
\label{eq:navier_stokes}
    \frac{\partial{\vect{u}}}{\partial t} + \nabla \cdot \vect{f(u)} - \nabla \cdot \vect{f^v(u,\nabla u)} = 0,
\end{equation}
where $\Omega \subset \mathcal{R}^d$ is the domain, $\vect{u}: \Omega \times [0\;T] \to \mathcal{R}^{d+2}$ is the conservative variable with the elements $\vect{u} = [\rho,\rho u_1, \dots,\rho u_d, \rho E]$, $\vect{f}$ is the convective flux and $\vect{f^v}$ is the viscous flux as defined below. $\rho$ is the density, $u_i$ is the velocity in the $i^{\text{th}}$ coordinate direction and $E$ is the specific total energy.  With $T$ as temperature, $p$ as pressure and $e$ as the internal energy density, we have $p=(\gamma-1)\rho e$, $e = c_v T$ and $E = e + \frac{1}{2} \sum_{i=1}^d u_i^2$, where $c_v$ is the specific heat at constant volume and $\gamma = 1.4$. 

The convective fluxes $\vect{f}^1,\dots,\vect{f}^d$ are given by
\begin{equation}
    \vect{f}^i(\vect{u}) = \begin{bmatrix}
        \rho u_i\\
        \rho u_i u_j + p\delta_{ij} \\
        u_i(\rho E + p)
    \end{bmatrix} \qquad j=1,\dots,d,
\end{equation}
while the dissipative fluxes are given by
\begin{equation}
    \vect{f}^{vi}(\vect{u},\vect{\nabla u}) = \begin{bmatrix}
        0\\
        \tau_{ji} \\
        \sum_{m=1}^d \tau_{mi}u_m - \kappa \frac{\partial T}{\partial x_i}
    \end{bmatrix} \qquad j=1,\dots,d.
\end{equation}
Here, $\kappa$ is the thermal conductivity and the viscous stress tensor $\tau$ is given by
\begin{equation}
    \tau_{ij} = \mu \left(\frac{\partial u_i}{\partial x_j} + \frac{\partial u_j}{\partial x_i}\right)
 +\delta_{ij}\lambda \sum_{m=1}^d \frac{\partial u_m}{\partial x_m}, \qquad 1\leq i,j\leq d. 
 \end{equation}
where $\mu$ is the viscosity and $\lambda = \frac{2}{3} \mu$ by Stokes hypothesis. We non-dimensionalize the Navier-Stokes equations with the freestream conditions \cite{masatsuka}.

The governing equation, \cref{eq:navier_stokes}, has an associated scalar convex entropy function, $U(\vect{u}) = -\rho s$, where $s = \log \left( \frac{p}{\rho^\gamma} \right)$ \cite{K_matrix_3D}. The entropy variables are $\vect{v}(\vect{u}) = U_{\vect{u}}$, and the mapping $\vect{u} \to \vect{v}$ is bijective and differentiable \cite{mock}. There exists a scalar flux functions $F^i(\vect{u})$, $i=1,\dots,d$, that satisfy $U_{\vect{u}} \vect{f^i_u} = F^i_{\vect{u}}$. Hence, $\vect{v}^T\frac{\partial \vect{f}^i}{\partial x_i} = \frac{\partial F^i}{\partial x_i}$ \cite{lax_entropy}. We also define the entropy potential, $\psi^i(\vect{u})$, that satisfies \cite{K_matrix_3D}  
\begin{equation}
    F^i(\vect{u}) = \vect{v}(\vect{u})^T\vect{f}^i(\vect{u})  - \psi^i(\vect{u}).
\end{equation}

The viscous flux, $\vect{f^v(u,\nabla u)}$, is linear in $\nabla u$ and can be written as $\vect{f^v(u,\nabla u)} = \vect{G(u)\nabla u}$. Using the differentability of the mapping $\vect{v} \to \vect{u}$, we can express it as $\vect{f^v(u,\nabla u)} = \vect{G(u)\frac{\partial u}{\partial v}\nabla v} = \vect{K}\vect{\nabla v}$, where $\vect{K} = \vect{G(u)\frac{\partial u}{\partial v}}$. It is known that $\vect{K}$ is symmetric positive semidefinite \cite{chan_entropystable_viscous}. Expressions for the tensor $\vect{K}$ can be found in \cite{K_matrix_3D,chan_entropystable_viscous}. 

Taking the inner product of \cref{eq:navier_stokes} with $\vect{v}(\vect{u})$ yields
\begin{equation}
\begin{aligned}
 &\int_\Omega \left(\frac{\partial U}{\partial t} + \nabla \cdot F - \vect{v}^T \nabla \cdot \vect{K\nabla v}\right) d\Omega = 0\\
 \implies & \int_\Omega \frac{\partial U}{\partial t} d\Omega + \int_{\partial \Omega} \left(F - \vect{v^TK\nabla v} \right) \cdot \vect{n} d\Gamma = -\int_\Omega \vect{\nabla v : K \nabla v} d\Omega \leq 0
\end{aligned}
\end{equation}

At the boundary $\partial \Omega$, we have $\vect{v^TK\nabla v} = \frac{1}{c_vT}\kappa \nabla T$ \cite{dalcin,chan_entropystable_viscous}. Therefore, the entropy inequality for the Navier-Stokes is
\begin{equation}   \label{eq:continuous_ns_entropy_inequality}
     \int_\Omega \frac{\partial U}{\partial t} d\Omega \leq \int_{\partial \Omega} \left(\frac{1}{c_vT}\kappa \nabla T - F \right) \cdot \vect{n} \;d\Gamma
\end{equation}

The goal is to have a discrete scheme that discretely satisfies the entropy inequality, \cref{eq:continuous_ns_entropy_inequality}. Entropy-stable treatment of the convective terms has been investigated extensively in the literature \cite{alex_nsfr,chan_entropystable_dg}. In this work, we focus on the entropy-stable discretizations of the viscous term, $-\nabla \cdot \vect{f}^v(\vect{u},\vect{\nabla u})$. 

\section{Notations}
As in the previous section, let $\Omega \subset \mathbb{R}^d$ be the domain of interest. We decompose the domain $\Omega$ into non-overlapping elements $\kappa$ and denote the set of these elements by $\mathcal{T}_h$. Let $\mathcal{E}_h$ be the set of all faces, $\mathcal{E}_{h,I}$ the set of all interior faces and $\mathcal{E}_{h,\partial}$ the set of all faces at the boundary $\partial \Omega$. Let $\hat{\kappa}$ be a standard reference element such that, for each physical element $\kappa$, there exists a diffeomorphism $\phi^\kappa : \kappa \to \hat{\kappa}$. Let the number of states in the governing equation be $n_s$. For the Navier-Stokes equations, we have $n_s = d+1$. For the purpose of discretization, we define the following finite-dimensional functional spaces:
\begin{equation}
    \mathcal{P}^p(\hat{\kappa}) = \{f:\hat{\kappa} \to \mathbb{R} : f \text{ is a polynomial of degree less than or equal to } p.\}
\end{equation}

\begin{equation}
     V_{h,p} = \{ w \in [L_2(\Omega)]^{n_s} : w|_k \circ \left(\phi^k\right)^{-1} \in [P^p(\hat{\kappa})]^{n_s} \quad \forall \kappa \in \mathcal{T}_h\}
\end{equation}

\begin{equation}
     \Sigma_{h,p} = \{ \tau \in [L_2(\Omega)]^{n_s \times d} : \tau|_k \circ \left(\phi^k\right)^{-1} \in [P^p(\hat{\kappa})]^{n_s \times d} \quad \forall \kappa \in \mathcal{T}_h\}
\end{equation}
where $w|_\kappa$ is the restriction of the function $w$ to element $\kappa$. We would like to stress that, for curvilinear grids, $V_{h,p}$ is not the space of polynomials on $\kappa$. Nevertheless, $V_{h,p}$ is still a finite-dimensional subspace of the Hilbert space $[L_2(\Omega)]^{n_s}$, and is hence a Hilbert space. The same holds for $\Sigma_{h,p}$.  

We define the average and jump operators at an edge $e$:
\begin{equation*}
    \{w\} = \begin{cases}
    & (w^+ + w^-)/2, \;\;  \text{in }\mathcal{E}_{h,I} \\
    & w^+, \;\; \text{on } \mathcal{E}_{h,\partial}
    \end{cases}  \qquad   [\![w]\!] = \begin{cases}
    & w^+n^+ + w^-n^-, \;\;  \text{in }\mathcal{E}_{h,I} \\
    & w^+n^+, \;\; \text{on } \mathcal{E}_{h,\partial}
    \end{cases}, \;\;  w \in \Pi_{\kappa \in \mathcal{T}_h} [L_2(\partial \kappa)]^{n_s}
\end{equation*}
and 
\begin{equation*}
    \{\tau\} = \begin{cases}
    & (\tau^+ + \tau^-)/2, \;\;  \text{in }\mathcal{E}_{h,I} \\
    & \tau^+, \;\; \text{on } \mathcal{E}_{h,\partial}
    \end{cases}  \qquad   [\![\tau]\!] = \begin{cases}
    & \tau^+ \cdot n^+ + \tau^- \cdot n^-, \;\;  \text{in }\mathcal{E}_{h,I} \\
    & \tau^+ \cdot n^+, \;\; \text{on } \mathcal{E}_{h,\partial}
    \end{cases}, \;\; \tau \in \Pi_{\kappa \in \mathcal{T}_h} [L_2(\partial \kappa)]^{n_s \times d}
\end{equation*}

With the above notations, the following identity is easy to verify by expanding the terms:
\begin{equation} 
\label{eq:avg_jump}
    \sum_{\kappa \in \mathcal{T}_h} \int_{\partial \kappa} w^T (\tau \cdot n) d\Gamma = \sum_{e \in \mathcal{E}_{h,I}} \int_e \{w\} \cdot [\![\tau]\!] d\Gamma + \sum_{e \in \mathcal{E}_{h}} \int_e \{\tau\} : [\![w]\!] d\Gamma
\end{equation}
Here, $n$ is the unit outward normal to an element face $\partial \kappa$.  

\subsection{Weight-adjusted projection and lifting operators}
\label{sec:weight_adjusted}
We denote $(u,v)_{\hat{\kappa}} = \int_{\hat{\kappa}} uv d\Omega$ to be the standard $L^2$ inner product on the reference element. Let $\Pi : L^2(\hat{\kappa}) \to \mathcal{P}^p(\hat{\kappa})$ be the usual projection operator on the reference element, defined by:
\begin{equation}
    (\Pi(u),\varphi)_{\hat{\kappa}} = (u,\varphi)_{\hat{\kappa}} \;\;\forall \varphi \in \mathcal{P}^p(\hat{\kappa}).
\end{equation}
If $J^k = \text{det}\left(\frac{\partial \left(\left(\phi^\kappa\right)^{-1}\right)}{\partial \xi} \right)$ is the Jacobian determinant of the mapping $\hat{\kappa} \to \kappa$, we have $(u,v)_{\kappa} = (J^ku,v)_{\hat{\kappa}}$ for any $u,v \in L^2(\kappa)$. 

For a valid mesh where $\exists \varepsilon_1,\varepsilon_2 > 0$ such that $\varepsilon_2>J^k > \varepsilon_1$ on $\Omega$, one can also equip the space $\mathcal{P}^p(\hat{\kappa})$ with the inner-product $<w_1,w_2> = \left(\frac{1}{J^k}w_1,w_2\right)$, with which $\mathcal{P}^p(\hat{\kappa})$ is still a Hilbert space. For a fixed $u \in L^2(\hat{\kappa})$, $(u,\varphi)_{\hat{\kappa}}$ is then a bounded linear functional $\mathcal{P}^p(\hat{\kappa}) \to \mathbb{R}$. Therefore, by the Riesz representation theorem, there exists a unique polynomial in $\mathcal{P}^p(\hat{\kappa})$, which we denote by $T^{-1}_{1/J^k}u$, that satisfies $<T^{-1}_{1/J^k}u,\varphi> = (u,\varphi)_{\hat{\kappa}}$, or
\begin{equation}
\label{eq:solving_for_T_inv_u}
    \left(\frac{1}{J^k} T^{-1}_{1/J^k} u,\varphi\right)_{\hat{\kappa}} = (u,\varphi)_{\hat{\kappa}}\;\;\forall \varphi \in \mathcal{P}^p(\hat{\kappa}).
\end{equation}
It has been proven that $T^{-1}_{1/J^k} u$ is a high-order approximation of $J^ku$ \cite{chan_weightadjusted_jcp} and is hence design order. If we restrict $T^{-1}_{1/J_k}$ to the finite-dimensional $\mathcal{P}^p(\hat{\kappa})$, then \cref{eq:solving_for_T_inv_u} shows that $T^{-1}_{1/J_k}: \mathcal{P}^p(\hat{\kappa}) \to \mathcal{P}^p(\hat{\kappa})$ is one-to-one and hence invertible. In fact, we have the following important properties of the operator $T^{-1}_{1/J_k}$ that we use in the subsequent sections \cite{chan_weightadjusted_math}:
\begin{enumerate}
    \item $T^{-1}_{1/J_k}$ is a bounded linear operator that is self-adjoint i.e. $\left(T^{-1}_{1/J_k}u,v\right)_{\hat{\kappa}} = \left(u,T^{-1}_{1/J_k}v\right)_{\hat{\kappa}} \;\;\forall u,v \in \mathcal{P}^p(\hat{\kappa})$.
    \item $T^{-1}_{1/J_k}$ is positive definite.
    \item Since $T^{-1}_{1/J_k}$ is invertible, we have $T^{-1}_{1/J_k}u=0$ only for $u=0\in \mathcal{P}^p(\hat{\kappa})$. Therefore, if we have $u,w \in \mathcal{P}^p(\hat{\kappa})$ that satsify $\left(T^{-1}_{1/J_k}u,\varphi \right)_{\hat{\kappa}} = \left(T^{-1}_{1/J_k}w,\varphi \right)_{\hat{\kappa}}\;\;\forall \varphi \in \mathcal{P}^p(\hat{\kappa})$, then $u = w$.
\end{enumerate}

The above properties imply that $\left(T^{-1}_{1/J_k}\cdot,\cdot\right)_{\hat{\kappa}}$ is an inner-product on $\mathcal{P}^p(\hat{\kappa})$, as also noted in \cite{chan_weightadjusted_math}. For conciseness, we denote it as $(\cdot,\cdot)_{T^{-1}_{1/J_k}} \equiv \left(T^{-1}_{1/J_k}\cdot,\cdot\right)_{\hat{\kappa}}$.

Next, we define the inner product on $\Sigma_{h,p}$ as $\sum_{\kappa \in \mathcal{T}_h}(u,v)_{T^{-1}_{1/J_k}} = \sum_{\kappa \in \mathcal{T}_h}\int_{\hat{\kappa}} T^{-1}_{1/J_k} u : v \;d\Omega$ for $u,v \in \Sigma_{h,p}$. Here, by "$:$" we mean $w:z = \sum_{i=1}^{n_s}\sum_{j=1}^d w_{ij}z_{ij}$. Equipped with this inner product, $\Sigma_{h,p}$ is a Hilbert space. 

We now define the following weight-adjusted projection and lift operators on an edge $e \in \mathcal{E}_h$:

\begin{enumerate}
 \item $r^e: [L_2(e)]^{n_s \times d} \to \Sigma_{h,p}$
    \begin{equation}
    \label{eq:lift_re}
        \sum_{\kappa \in \mathcal{T}_h} \left( r^e(\phi) : \tau^h \right)_{T^{-1}_{1/J_k}} = - \int_e \phi : \{\tau^h\} d\Gamma \quad \forall \tau^h \in \Sigma_{h,p}
    \end{equation}
    Note: For a fixed $\phi \in [L_2(e)]^{n_s \times d}$, the right-hand side is a linear functional in $\tau^h$. With the inner product on $\Sigma_{h,p}$ defined above and the Riesz representation theorem, existence of $r^e(\phi) \in \Sigma_{h,p}$ follows.
    \item $l^e: [L_2(e)]^{n_s} \to \Sigma_{h,p}$
    \begin{equation}
    \label{eq:lift_le}
        \sum_{\kappa \in \mathcal{T}_h} \left( l^e(q) : \tau^h \right)_{T^{-1}_{1/J_k}} = - \int_e q^T [\![\tau^h]\!] d\Gamma \quad \forall \tau^h \in \Sigma_{h,p}
    \end{equation}
    \item $\mathcal{G}^h: [L_2(\Omega)]^{m \times d} \to \Sigma_{h,p}$
\begin{equation}
\label{eq:Gh}
    \sum_{\kappa \in \mathcal{T}_h} \left(\mathcal{G}^h(\phi) : \tau^h \right)_{T^{-1}_{1/J_k}} = \sum_{\kappa \in \mathcal{T}_h} \int_{\kappa} \phi : \tau^h d\Omega \quad \forall \tau^h \in \Sigma_{h,p}
\end{equation}
\item $r = \sum_{e\in \mathcal{E}_h}r^e$, $\;\;r_\mathcal{I} = \sum_{e\in \mathcal{E}_{h,I}}r^e$ and $r_\Gamma = \sum_{e\in \mathcal{E}_{h,\partial}}r^e$. Similar notations are used for $l^e$.
\end{enumerate} 

\section{Entropy-stable viscous discretizations}
Based on the discussion in \Cref{sec:entropy_stable_ns}, the Navier-Stokes equations can be written as 
\begin{subequations}
    \begin{align}
    & \theta = \nabla v, \label{eq:theta}\\
    & \sigma = K\theta, \label{eq:sigma} \\
    & \frac{\partial u}{\partial t} + \nabla \cdot f(u) - \nabla \cdot \sigma = 0. \label{eq:ns}
    \end{align}
\end{subequations}

We discretize \cref{eq:theta} as: find $\theta^h \in \Sigma_{h,p}$ that satisfies:
\begin{equation}
    \label{eq:theta_h} 
    \sum_{\kappa \in \mathcal{T}_h} \left(\theta^h, \tau^h\right)_{T^{-1}_{1/J_k}} = \sum_{\kappa \in \mathcal{T}_h} \int_{\partial \kappa} (\tau^h\cdot n)\cdot (\hat{v}-v^h) d\Gamma + \int_\kappa \nabla v^h : \tau^h d\Omega \quad \forall \tau^h \in \Sigma_{h,p},
\end{equation}
where $\hat{v}$ is the numerical flux for the entropy variable, which will be specific to the discretization scheme, and $v^h$ is the weight-adjusted projection of the entropy variable, defined as $\left(v^h,\varphi \right)_{T^{-1}_{1/J_k}} = \int_\kappa \vect{v(u)}^T\varphi\; d\Omega\;\; \forall \varphi \in V_{h,p}$.

Similarly, we approximate \cref{eq:sigma} as: find $\sigma_h \in \Sigma_{h,p}$ that satisifes:
\begin{equation}
\label{eq:sigma_h}
    \sum_{\kappa \in \mathcal{T}_h} \left(\sigma^h, \zeta^h \right)_{T^{-1}_{1/J_k}} =  \sum_{\kappa \in \mathcal{T}_h} \int_\kappa K(v^h)\theta^h : \zeta^h d\Omega  \quad \forall \zeta^h \in \Sigma_{h,p}.
\end{equation}
$\theta^h$ and $\sigma^h$ exist from the Reisz representation theorem if $\hat{v} \in [L_2(e)]^{n_s}$ and $K\theta^h \in [L_2(\Omega)]^{n_s \times d}$.

Finally, we discretize \cref{eq:ns} as: find $u^h \in V_{h,p}$ that satisfies:
\begin{equation}
    \sum_{\kappa \in \mathcal{T}_h}\left(\frac{\partial u^h}{\partial t},w^h\right)_{T^{-1}_{1/J_k}} + {B}^c(u^h,w^h) + B^v(u^h,w^h) = 0 \quad \forall w^h \in V_{h,p}
\end{equation}
where $B^c(u^h,w^h)$ is an entropy-stable discretization of the convective term $\left(\nabla \cdot f(u^h),w^h\right)_\Omega$ and 
\begin{equation}
\label{eq:Bv_original}
    B^v(u^h,w^h) = \sum_{\kappa \in \mathcal{T}_h} \left( -\int_{\partial \kappa} (\hat{\sigma}\cdot n)w^h d\Gamma + \int_\kappa \nabla w^h: \sigma^h \;d\Omega \right). 
\end{equation}
We note an important difference between the conventional DG formulation and the weight-adjusted formulation: conventional DG has $\left(\frac{\partial u^h}{\partial t},w^h\right)_\kappa$ in place of $\left(\frac{\partial u^h}{\partial t},w^h\right)_{T^{-1}_{1/J_k}}$. However, as noted in \Cref{sec:weight_adjusted}, $\left(u^h,w^h\right)_{T^{-1}_{1/J_k}}$ is a high-order approximation of $\left(J^ku^h,w^h\right)_{\hat{\kappa}} = \left(u^h,w^h\right)_\kappa$ \cite{chan_weightadjusted_jcp}. 

Using the identity in \cref{eq:avg_jump} to the right-hand side of \cref{eq:theta_h}, and the fact that any valid numerical flux satisfies $\{\hat{v}\} = \hat{v}$ on $\mathcal{E}_h$ and $[\![\hat{v}]\!] = 0$ in $\mathcal{E}_{h,I}$, we obtain
\begin{equation}
\begin{split}
    \sum_{\kappa \in \mathcal{T}_h} \left(\theta^h, \tau^h\right)_{T^{-1}_{1/J_k}}  =& \sum_{\kappa \in \mathcal{T}_h} \int_\kappa \nabla v^h : \tau^h d\Omega + \sum_{e \in \mathcal{E}_{h,I}} \int_e (\hat{v} - \{v^h\})\cdot [\![\tau^h]\!] d\Gamma \\
    &- \sum_{e \in \mathcal{E}_h} \int_e \{\tau^h\} \cdot [\![v^h]\!] d\Gamma + \sum_{e \in \mathcal{E}_{h,\partial}} \int_e \{\tau^h\} \cdot [\![\hat{v}]\!] d\Gamma
    \end{split}
\end{equation}

Using the weight-adjusted lift and projection operators defined in \Cref{sec:weight_adjusted}, the above equation can be written as 

\begin{equation}
\begin{split}
    \sum_{\kappa \in \mathcal{T}_h} \left(\theta^h, \tau^h\right)_{T^{-1}_{1/J_k}}  = &\sum_{\kappa \in \mathcal{T}_h}  \left(\mathcal{G}^h(\nabla v^h), \tau^h\right)_{T^{-1}_{1/J_k}} + \sum_{e\in\mathcal{E}_{h,I}} \sum_{\kappa \in \mathcal{T}_h}  \left(l^e(\{v^h\} - \hat{v}), \tau^h\right)_{T^{-1}_{1/J_k}} \\
    & + \sum_{e\in\mathcal{E}_{h}} \sum_{\kappa \in \mathcal{T}_h} \left(r^e([\![v^h]\!], \tau^h\right)_{T^{-1}_{1/J_k}} -  \sum_{e\in\mathcal{E}_{h,\partial}}\sum_{\kappa \in \mathcal{T}_h}  \left(r^e(\hat{v}\cdot n), \tau^h\right)_{T^{-1}_{1/J_k}} \;\; \forall \tau^h \in \Sigma_{h,p}.
    \end{split}
\end{equation}
Therefore, as noted in \Cref{sec:weight_adjusted} about the properties of $(\cdot,\cdot)_{T^{-1}_{1/J_k}}$, we have the equality
\begin{equation}
\begin{split}
\label{eq:final_theta_h}
    \theta^h &= \mathcal{G}^h(\nabla v^h) + l_{\mathcal{I}}(\{v^h\} - \hat{v}) + r([\![v^h]\!]) - r_\Gamma(\hat{v} n) \\
    &= \mathcal{G}^h(\nabla v^h) + l_{\mathcal{I}}(\{v^h\} - \hat{v}) + r_{\mathcal{I}}([\![v^h]\!]) - r_\Gamma((\hat{v}-v^h) n).
    \end{split}
\end{equation}

From \cref{eq:sigma_h}, we have
\begin{equation}
\label{eq:sigma_h_projection}
    \sigma_h = \mathcal{G}^h\left(K\theta^h\right).
\end{equation}

Finally, using the definition of $\mathcal{G}^h$ from \cref{eq:Gh}, we can write $\mathcal{B}^v$ in \cref{eq:Bv_original} as
\begin{equation}
\begin{split}
\label{eq:Bv_2}
        B^v(u^h,w^h)&=\sum_{\kappa \in \mathcal{T}_h} \left( -\int_{\partial \kappa} (\hat{\sigma}\cdot n)w^h d\Gamma + \int_\kappa \nabla w^h: \sigma^h \;d\Omega \right)  \\
        &=\sum_{\kappa \in \mathcal{T}_h} \left( -\int_{\partial \kappa} (\hat{\sigma}\cdot n)w^h d\Gamma + \left(\mathcal{G}^h(\nabla w^h), \sigma^h\right)_{T^{-1}_{1/J_k}} \right).
\end{split}
\end{equation}
Applying \cref{eq:avg_jump} to the first term in the above equation, while noting that $\{\hat{\sigma}\} = \hat{\sigma}$ on $\mathcal{E}_h$ and $[\![\hat{\sigma}]\!] = 0$ on $\mathcal{E}_{h,I}$, yields
\begin{equation}
    \begin{split}
    \label{eq:Bv3}
        B^v(u^h,w^h)= \sum_{\kappa \in \mathcal{T}_h} \left(\mathcal{G}^h(\nabla w^h), \sigma^h\right)_{T^{-1}_{1/J_k}} - \sum_{e \in \mathcal{E}_h} \int_e \hat{\sigma} : [\![w^h]\!] \; d\Gamma
    \end{split}
\end{equation}
For the first term above, using \cref{eq:sigma_h_projection}, \cref{eq:Gh} and the fact that $K(v^h)$ is a symmetric tensor yields 
\begin{equation}
\begin{split}
\label{eq:Bv_interim}
  \sum_{\kappa \in \mathcal{T}^h} \left(\mathcal{G}^h(\nabla w^h), \sigma^h\right)_{T^{-1}_{1/J_k}} &= \sum_{\kappa \in \mathcal{T}^h}\left(\mathcal{G}^h(\nabla w^h), \mathcal{G}^h(K\theta^h)\right)_{T^{-1}_{1/J_k}} = \sum_{\kappa \in \mathcal{T}^h} \int_\kappa \mathcal{G}^h(\nabla w^h) : K\theta^h d\Omega \\&= \sum_{\kappa \in \mathcal{T}^h} \int_\kappa K\mathcal{G}^h(\nabla w^h) : \theta^h d\Omega.
  \end{split}
\end{equation}
Therefore, substituting \cref{eq:Bv_interim} into \cref{eq:Bv3} and using \cref{eq:final_theta_h} yields
\begin{equation}
    \begin{split}
    \label{eq:Bv4}
        B^v(u^h,w^h)&= \sum_{\kappa \in \mathcal{T}^h} \int_\kappa K\mathcal{G}^h(\nabla w^h) : \theta^h d\Omega - \sum_{e \in \mathcal{E}_h} \int_e \hat{\sigma} : [\![w^h]\!] \; d\Gamma \\
        &= \sum_{\kappa \in \mathcal{T}_h} \int_\kappa K\mathcal{G}^h(\nabla w^h): \left( \mathcal{G}^h(\nabla v^h) + l_{\mathcal{I}}(\{v^h\}- \hat{v}) + r_{\mathcal{I}}([\![v^h]\!]) - r_\Gamma((\hat{v}-v^h) n)\right) d\Omega \\&\;\;\;- \sum_{e \in \mathcal{E}_h} \int_e \hat{\sigma} : [\![w^h]\!] \; d\Gamma.
    \end{split}
\end{equation}

For the ease of notation, we split \cref{eq:Bv4} into two components,
\begin{equation}
     B^v(u^h,w^h) =  B_{\mathcal{I}}^v(u^h,w^h) +  B_{\Gamma}^v(u^h,w^h),
\end{equation}
where
\begin{equation}
\begin{split}
\label{eq:BI_final}
    B_{\mathcal{I}}^v(u^h,w^h) =& \sum_{\kappa \in \mathcal{T}_h} \int_\kappa K\mathcal{G}^h(\nabla w^h): \mathcal{G}^h(\nabla v^h) d\Omega \\
    &+ \sum_{\kappa \in \mathcal{T}_h}\int_\kappa K\mathcal{G}^h(\nabla w^h): \left(r_{\mathcal{I}}([\![v^h]\!])+l_{\mathcal{I}}(\{v^h\}- \hat{v})\right) d\Omega \\
    &-\sum_{e \in \mathcal{E}_{h,I}} \int_e \hat{\sigma} : [\![w^h]\!] \; d\Gamma,
\end{split}
\end{equation}
and
\begin{equation}
\label{eq:B_gamma_final}
      B_{\Gamma}^v(u^h,w^h) = -\sum_{e \in \mathcal{E}_{h,\partial}}\sum_{\kappa \in \mathcal{T}_h}\int_\kappa K\mathcal{G}^h(\nabla w^h) : r^e\left((\hat{v} - v^h) \otimes n\right) \;d\Omega -\sum_{e \in \mathcal{E}_{h,\partial}} \int_e \hat{\sigma} : [\![w^h]\!] \; d\Gamma.  
\end{equation}
In \cref{eq:B_gamma_final}, we have substituted in $r_\Gamma = \sum_{e \in \mathcal{E}_{h,\partial}} r^e$.
\begin{remark}
    We can remove the weight-adjusted projections from the terms involving the test function $w^h$ by noting that, for any $a^h \in \Sigma_{h,p}$, we have 
    \begin{equation*}
        \int_\kappa K\mathcal{G}^h(\nabla w^h) : a^h d\Omega = \int_\kappa \mathcal{G}^h(\nabla w^h) : Ka^h d\Omega =\left(\mathcal{G}^h(\nabla w^h),\mathcal{G}^h\left(Ka^h\right)\right)_{T^{-1}_{`/J_k}} = \int_\kappa \nabla w^h : \mathcal{G}^h\left(Ka^h\right) d\Omega.
    \end{equation*}
\end{remark}
\begin{remark}
    The derivations and proofs that follow use only the equalities in \cref{eq:lift_re,eq:lift_le,eq:Gh}. Thus, even if the integration is inexact, the subsequent proofs hold as long as \cref{eq:lift_re,eq:lift_le,eq:Gh} are satisfied with the inexact integration. 
\end{remark}

\begin{remark}
    The weight-adjusted formulation is known to be of design order, as mathematically proven and numerically demonstrated in the literature  \cite{chan_weightadjusted_jcp,chan_weightadjusted_math,alex_nsfr}. Therefore, the above weight-adjusted formulation for the viscous discretization is expected to preserve the high-order convergence of the DG scheme.
\end{remark}

\subsection{Weight-adjusted entropy-stable BR2 (BR2-E)}
For the modified BR2 formulation, which we refer to as BR2-E, we use the following specific numerical fluxes:

\begin{equation}
\label{eq:bc_br2}
    \hat{v} = \begin{cases}
        &\{v^h\},\quad \text{in }\mathcal{E}_{h,I} \\
        & 2v^h - v_\Gamma(v^h),\quad \text{on }\mathcal{E}_{h,\partial}
    \end{cases} \qquad \hat{\sigma} = \begin{cases}
        &\{\mathcal{G}^h\left(K \mathcal{G}^h(\nabla v^h) + \eta_e K r^e([\![v^h]\!])\right)  \},\quad \text{in }\mathcal{E}_{h,I} \\
        & \sigma_\Gamma(\mathcal{G}^h(K\mathcal{G}^h(\nabla v^h))),\quad \text{on }\mathcal{E}_{h,\partial}
    \end{cases}
\end{equation}
where $\eta_e>0$ is a penalization. Note that the fluxes in the interior faces are similar to the usual BR2 formulation \cite{bassi_rebay2,May_viscous_entropystableDG}, with the difference of using a weight-adjusted projection on the gradient $\nabla v^h$, along with the weight-adjusted lifting operators. Also note that the fluxes at the boundary are modified. In particular, we note that the boundary flux $\hat{v} = 2v^h - v_\Gamma(v^h)$ leads to a consistent discretization if $v_\Gamma$ is the boundary condition for the entropy variable. The same holds for $\hat{\sigma}$ if $\sigma_\Gamma$ is consistent with the boundary condition for the viscous flux.  

Substituting \cref{eq:bc_br2} into \cref{eq:BI_final}, and noting that $l_{\mathcal{I}}(\{v^h\} - \hat{v}) = 0$, we get

\begin{equation}
\begin{split}
\label{eq:BI_br2_1}
    B_{\mathcal{I}}^{v,\text{BR2-E}}(u^h,w^h) =&\sum_{\kappa \in \mathcal{T}_h} \int_\kappa K\mathcal{G}^h(\nabla w^h): \mathcal{G}^h(\nabla v^h) d\Omega  + \sum_{\kappa \in \mathcal{T}_h} \int_{\kappa} K\mathcal{G}^h(\nabla w^h) : r_{\mathcal{I}}([\![v^h]\!]) d\Omega \\ & - \sum_{e \in \mathcal{E}_{h,I}} \int_e \hat{\sigma} : [\![w^h]\!] \; d\Gamma .
\end{split}
\end{equation}
The last term becomes
\begin{equation}
\begin{split}
\label{eq:sigmahat_br2}
    - \sum_{e \in \mathcal{E}_{h,I}} \int_e \hat{\sigma} : [\![w^h]\!] \; d\Gamma &= - \sum_{e \in \mathcal{E}_{h,I}} \int_e \{\mathcal{G}^h(K \mathcal{G}^h(\nabla v^h) + \eta_e K r^e([\![v^h]\!]))  \} : [\![w^h]\!] \; d\Gamma\\
    &= \sum_{e \in \mathcal{E}_{h,I}}\sum_{\kappa \in \mathcal{T}_h} \left(r^e([\![w^h]\!]),\;\mathcal{G}^h\left(K \mathcal{G}^h(\nabla v^h) + \eta_e K r^e([\![v^h]\!])\right) \right)_{T^{-1}_{1/{J_k}}}
     \\
    &= \sum_{\kappa \in \mathcal{T}_h} \int_{\kappa} K \mathcal{G}^h(\nabla v^h) : r_{\mathcal{I}}([\![w^h]\!]) d\Omega + \sum_{e \in \mathcal{E}_{h,I}} \sum_{\kappa \in \mathcal{T}_h} \int_\kappa \eta_e K r^e([\![v^h]\!]) : r^e([\![w^h]\!]) d\Omega 
    \end{split}
\end{equation}
where we have used \cref{eq:lift_re} and \cref{eq:Gh}. Substituting \cref{eq:sigmahat_br2} into \cref{eq:BI_br2_1} yields
\begin{equation}
    \begin{split}
\label{eq:BI_br2_final}
    B_{\mathcal{I}}^{v,\text{BR2-E}}(u^h,w^h)
    =&\sum_{\kappa \in \mathcal{T}_h} \int_\kappa K\mathcal{G}^h(\nabla w^h): \mathcal{G}^h(\nabla v^h) d\Omega  + \sum_{\kappa \in \mathcal{T}_h} \int_{\kappa} K\mathcal{G}^h(\nabla w^h) : r_{\mathcal{I}}([\![v^h]\!]) d\Omega \\
    &+\sum_{\kappa \in \mathcal{T}_h} \int_{\kappa} K \mathcal{G}^h(\nabla v^h) : r_{\mathcal{I}}([\![w^h]\!]) d\Omega + \sum_{e \in \mathcal{E}_{h,I}} \sum_{\kappa \in \mathcal{T}_h} \int_\kappa \eta_e K r^e([\![v^h]\!]) : r^e([\![w^h]\!]) d\Omega.
    \end{split}
\end{equation}

On the other hand, from \cref{eq:B_gamma_final} and \cref{eq:bc_br2}, we have
\begin{equation}
\label{eq:B_gamma_br2}
      B_{\Gamma}^{v,\text{BR2-E}}(u^h,w^h) =  \sum_{e \in \mathcal{E}_{h,\partial}}\sum_{\kappa \in \mathcal{T}_h}\int_\kappa K\mathcal{G}^h(\nabla w^h) : r^e\left((v_\Gamma(v^h) - v^h) \otimes n\right) \;d\Omega\;   -\sum_{e \in \mathcal{E}_{h,\partial}} \int_e \sigma_\Gamma(\mathcal{G}^h(K\mathcal{G}^h(\nabla v^h))) : [\![w^h]\!] \; d\Gamma.  
\end{equation}

The next theorem shows that $\mathcal{B}^{v,\text{BR2-E}}_\mathcal{I}(u^h,v^h) \geq 0$. The proof is similar to the energy stability proof of BR2 \cite{unified_viscous_dg}.
\begin{theorem}
\label{thm:BI_br2}
    Let $N_f$ be the number of faces of an element. If $\eta_e>N_f$, then $\mathcal{B}^{v,\text{BR2-E}}_\mathcal{I}(u^h,v^h) \geq 0$.
\end{theorem}
\begin{proof}
   We have
   \begin{equation}
     \begin{split}
\label{eq:BI_br2_vv}
    B_{\mathcal{I}}^{v,\text{BR2-E}}(u^h,v^h)
    =&\sum_{\kappa \in \mathcal{T}_h} \int_\kappa K\mathcal{G}^h(\nabla v^h): \mathcal{G}^h(\nabla v^h) d\Omega  + 2\sum_{\kappa \in \mathcal{T}_h} \int_{\kappa} K\mathcal{G}^h(\nabla v^h) : r_{\mathcal{I}}([\![v^h]\!]) d\Omega \\
    & + \sum_{e \in \mathcal{E}_{h,I}} \sum_{\kappa \in \mathcal{T}_h} \int_\kappa \eta_e K r^e([\![v^h]\!]) : r^e([\![v^h]\!]) d\Omega.
    \end{split}
   \end{equation}
\Cref{eq:BI_br2_vv} is similar to \cite[Eq. 54]{May_viscous_entropystableDG}, with three important differences: First, the lifting operators $r^e$ and $r_{\mathcal{I}}$ in the current work are defined using the weight-adjusted formulation, \cref{eq:lift_re}. Secondly, the second term in \cref{eq:BI_br2_vv} has $r_{\mathcal{I}}$, and the third term has a summation over the interior edges $\mathcal{E}_{h,I}$, whereas in \cite[Eq. 54]{May_viscous_entropystableDG} the second term has $r$ and the third term has summation over all the edges $\mathcal{E}_h$ and lastly, we have the weight-adjusted projection $\mathcal{G}^h(\nabla v^h)$ instead of $\nabla v^h$. Nevertheless, since \cite[Eqs. 54-56]{May_viscous_entropystableDG} is agnostic to the definition of the lifting operators, we can apply these equations directly with the summation over $\mathcal{E}_{h,I}$ instead, proving $B_{\mathcal{I}}^{v,\text{BR2-E}}(u^h,v^h)\geq0$ if $\eta_e > N_f$. 
\end{proof}

\begin{theorem}
\label{thm:B_gamma_br2}
   Denote $\tilde{\sigma}^h = \mathcal{G}^h(K\mathcal{G}^h(\nabla v^h))$. Let
    \begin{equation}
        v_\Gamma(v^h) = \frac{v^+(v^h)+v^h}{2}
    \end{equation}
    and 
    \begin{equation}
        \sigma_\Gamma(\tilde{\sigma}^h) = \frac{\sigma^+(\tilde{\sigma}^h) + \tilde{\sigma}^h}{2}.
    \end{equation}
    If the functions $v^+(\cdot)$ and $\sigma^+(\cdot)$ are chosen as specified in \cite{chan_entropystable_viscous}, then $\mathcal{B}^{v,\text{BR2-E}}_\Gamma(u^h,v^h)$ is a consistent discretization of the entropy boundary integral in \cref{eq:continuous_ns_entropy_inequality}, for all the boundary conditions presented in \cite{chan_entropystable_viscous}.
\end{theorem}
\begin{proof}
    From \cref{eq:B_gamma_br2}, we have
    \begin{equation}
    \begin{split}
      B_{\Gamma}^{v,\text{BR2-E}}(u^h,v^h) &=  \sum_{e \in \mathcal{E}_{h,\partial}}\sum_{\kappa \in \mathcal{T}_h}\int_\kappa K\mathcal{G}^h(\nabla v^h) : r^e\left((v_\Gamma(v^h) - v^h) \otimes n\right) \;d\Omega\;   -\sum_{e \in \mathcal{E}_{h,\partial}} \int_e \sigma_\Gamma(\mathcal{G}^h(K\mathcal{G}^h(\nabla v^h))) : [\![v^h]\!] \; d\Gamma \\
       &=  \sum_{e \in \mathcal{E}_{h,\partial}}\sum_{\kappa \in \mathcal{T}_h}\left(\mathcal{G}^h\left(K\mathcal{G}^h(\nabla v^h)\right),r^e\left((v_\Gamma(v^h) - v^h) \otimes n\right) \right)_{T^{-1}_{1/J_k}}   -\sum_{e \in \mathcal{E}_{h,\partial}} \int_e \sigma_\Gamma(\mathcal{G}^h(K\mathcal{G}^h(\nabla v^h))) : [\![v^h]\!] \; d\Gamma \\
      &=-\sum_{e \in \mathcal{E}_{h,\partial}} \int_e \mathcal{G}^h(K\mathcal{G}^h(\nabla v^h)) : (v_\Gamma(v^h)-v^h) \otimes n \;d\Gamma -\sum_{e \in \mathcal{E}_{h,\partial}} \int_e \sigma_\Gamma(\mathcal{G}^h(K\mathcal{G}^h(\nabla v^h))) : [\![v^h]\!] \; d\Gamma \\
      &= -\sum_{e \in \mathcal{E}_{h,\partial}} \int_e \tilde{\sigma}^h : (v_\Gamma(v^h)-v^h) \otimes n \;d\Gamma -\sum_{e \in \mathcal{E}_{h,\partial}} \int_e \sigma_\Gamma(\tilde{\sigma}^h) : [\![v^h]\!] \; d\Gamma 
      \end{split}
\end{equation}
Substituting $v_\Gamma$ and $\sigma_\Gamma$ yields
\begin{equation}
\label{eq:B_gamma_JC}
     B_{\Gamma}^{v,\text{BR2-E}}(u^h,v^h) = -\sum_{e \in \mathcal{E}_{h,\partial}} \frac{1}{2}\int_e \tilde{\sigma}^h : (v^+(v^h)-v^h) \otimes n \;d\Gamma -\sum_{e \in \mathcal{E}_{h,\partial}} \int_e \frac{\sigma^+(\tilde{\sigma}^h) + \tilde{\sigma}^h}{2} : v^h \otimes n\; d\Gamma.
\end{equation}
    We note that \cref{eq:B_gamma_JC} is the same as the boundary integral in \cite[Eq. 22]{chan_entropystable_viscous}, moved to the left-hand side. Therefore, proofs in \cite{chan_entropystable_viscous} apply directly, showing that $\mathcal{B}^{v,\text{BR2-E}}_\Gamma(u^h,v^h)$ is a consistent discretization of the entropy boundary integral.  
\end{proof}

The above two proofs show that the weight-adjusted BR2-E formulation is entropy-stable. Formally,
\begin{theorem}
\label{thm:entropy_stability_br2}
    Consider the BR2-E discretization scheme with the penalization $\eta_e > N_f$:
    \begin{equation}
    \sum_{\kappa \in \mathcal{T}_h}\left(\frac{\partial u^h}{\partial t},w^h\right)_{T^{-1}_{1/J_k}} + {B}^c(u^h,w^h) + B^{v,\text{BR2-E}}(u^h,w^h) = 0 \quad \forall w^h \in V_{h,p},
\end{equation}
where $B^{v,\text{BR2-E}}(u^h,w^h) = B_{\mathcal{I}}^{v,\text{BR2-E}}(u^h,w^h) + B_\Gamma^{v,\text{BR2-E}}(u^h,w^h)$, with $B_{\mathcal{I}}^{v,\text{BR2-E}}(u^h,w^h)$ given by \cref{eq:BI_br2_final} and $B_\Gamma^{v,\text{BR2-E}}(u^h,w^h)$ given by \cref{eq:B_gamma_br2}. If $\mathcal{B}^c(u^h,w^h)$ is an entropy-stable discretization of the convective term, the scheme is entropy-stable. 
\end{theorem}
\begin{proof}
    By the construction of $v^h$, we have $\left(v^h,\varphi \right)_{T^{-1}_{1/J_k}} = \int_\kappa \vect{v(u)}^T\varphi\; d\Omega\;\; \forall \varphi \in V_{h,p}$ \cite{chan_weightadjusted_jcp}. Therefore, substituting $v^h$ in place of $w^h$ yields
    \begin{equation}
    \label{eq:entropy_discrete_inequality}
        \begin{split}
            &\sum_{\kappa \in \mathcal{T}_h}\left(\frac{\partial u^h}{\partial t},v^h\right)_{T^{-1}_{1/J_k}} = - {B}^c(u^h,v^h) - B^{v,\text{BR2-E}}_{\mathcal{I}}(u^h,v^h)- B^{v,\text{BR2-E}}_\Gamma(u^h,v^h)\\
            &\sum_{\kappa \in \mathcal{T}_h} \int_\kappa \vect{v}(u^h)^T\frac{\partial u^h}{\partial t} d\Omega \leq - B^{v,\text{BR2-E}}_\Gamma(u^h,v^h)\\
            &\sum_{\kappa \in \mathcal{T}_h} \int_\kappa \frac{\partial U(u^h)}{\partial t} d\Omega \leq - B^{v,\text{BR2-E}}_\Gamma(u^h,v^h)
        \end{split}
    \end{equation}
    since $- B^{v,\text{BR2-E}}_{\mathcal{I}}(u^h,v^h)\leq 0$ from \Cref{thm:BI_br2} and $- B^c(u^h,v^h)\leq 0$ because $\mathcal{B}^c$ is entropy-stable. From \Cref{thm:B_gamma_br2}, we can conclude that the last inequality in \cref{eq:entropy_discrete_inequality} is a consistent discretization of the continuous entropy inequality, \cref{eq:continuous_ns_entropy_inequality}, for all the boundary conditions presented in \cite{chan_entropystable_viscous}.
\end{proof}

\subsection{Weight-adjusted entropy-stable CDG (CDG-E)}
The compact discontinuous Galerkin (CDG) is a modification of the LDG scheme in order to obtain a compact stencil. It was introduced in \cite{CDG_scheme}. While the original formulation has numerical fluxes in terms of the primal solution $u^h$, we modify it to be applied to the entropy variable $v^h$. The modified viscous numerical fluxes for our purpose of entropy stability for CDG-E are
\begin{equation}
\label{eq:bc_cdg}
    \hat{v} = \begin{cases}
        &\{v^h\} - \vect{C_{12}}\cdot [\![v^h]\!],\quad \text{in }\mathcal{E}_{h,I} \\
        & 2v^h - v_\Gamma(v^h),\quad \text{on }\mathcal{E}_{h,\partial}
    \end{cases} \qquad \hat{\sigma} = \begin{cases}
        &\{\sigma^e_h\} - C_{11}[\![v^h]\!] + \vect{C_{12}}[\![\sigma^e_h]\!],\quad \text{in }\mathcal{E}_{h,I} \\
        & \sigma_\Gamma(\mathcal{G}^h(K\mathcal{G}^h(\nabla v^h))),\quad \text{on }\mathcal{E}_{h,\partial}
    \end{cases}
\end{equation}
where 
\begin{equation}
    \sigma^e_h = \mathcal{G}^h\left(K\mathcal{G}^h(\nabla v^h)\right) + \overline{\sigma}^e_h
\end{equation}
and
\begin{equation}
    \overline{\sigma}^e_h = \mathcal{G}^h\left(Kr^e([\![v^h]\!]) + Kl^e(\vect{C_{12}}\cdot [\![v^h]\!]) \right).
\end{equation}

Here, $C_{11}>0$ is a scalar constant and $\vect{C_{12}} = \frac{1}{2}\left(S_{\kappa^-}^{\kappa^+}n^+ + S_{\kappa^+}^{\kappa^-}n^-  \right)$, with the switch $S_{\kappa^+}^{\kappa^-} \in \{0,1\}$ on each element face. The switch satisfies $S_{\kappa^+}^{\kappa^-} + S_{\kappa^-}^{\kappa^+}=1$. For instance, one can use the global numbering of cells and set $\vect{C_{12}} = \frac{1}{2}n^+$ at the face between $\kappa^+$ and $\kappa^-$, if the cell id of $\kappa^+$ is less than the cell id of $\kappa^-$. 

Substituting \cref{eq:bc_cdg} into \cref{eq:BI_final} yields
\begin{equation}
\begin{split}
\label{eq:BI_cdg_1}
    B_{\mathcal{I}}^{v,\text{CDG-E}}(u^h,w^h) =&\sum_{\kappa \in \mathcal{T}_h} \int_\kappa K\mathcal{G}^h(\nabla w^h): \mathcal{G}^h(\nabla v^h) d\Omega + \sum_{\kappa \in \mathcal{T}_h} \int_\kappa  K\mathcal{G}^h(\nabla w^h):\left(r_{\mathcal{I}}([\![v^h]\!]) + l_{\mathcal{I}}(\vect{C_{12}}\cdot [\![v^h]\!])\right) d\Omega \\
    & - \sum_{e \in \mathcal{E}_{h,I}} \int_e \hat{\sigma} : [\![w^h]\!] \; d\Gamma,
\end{split}
\end{equation}
where
\begin{equation}
    \begin{split}
         - \sum_{e \in \mathcal{E}_{h,I}} \int_e \hat{\sigma} : [\![w^h]\!] \; d\Gamma =&  - \sum_{e \in \mathcal{E}_{h,I}} \int_e \{\sigma^e_h\} : [\![w^h]\!] \; d\Gamma +  \sum_{e \in \mathcal{E}_{h,I}} \int_e C_{11}[\![v^h]\!] : [\![w^h]\!] \; d\Gamma \\
         &- \sum_{e \in \mathcal{E}_{h,I}} \int_e \left(\vect{C_{12}}\cdot  [\![w^h]\!]\right) \cdot [\![\sigma^e_h]\!] \; d\Gamma \\
         =& \sum_{e \in \mathcal{E}_{h,I}}\sum_{\kappa \in \mathcal{T}_h} \left(r^e([\![w^h]\!]),\sigma^e_h\right)_{T^{-1}_{1/J_k}} +  \sum_{e \in \mathcal{E}_{h,I}} \int_e C_{11}[\![v^h]\!] : [\![w^h]\!] \; d\Gamma \\
         &+ \sum_{e \in \mathcal{E}_{h,I}}\sum_{\kappa \in \mathcal{T}_h} \left(l^e\left(\vect{C_{12}}\cdot[\![w^h]\!]\right),\sigma^e_h\right)_{T^{-1}_{1/J_k}} \\
         = & \sum_{e \in \mathcal{E}_{h,I}}\sum_{\kappa \in \mathcal{T}_h} \left(r^e([\![w^h]\!]) + l^e\left(\vect{C_{12}}\cdot[\![w^h]\!]\right),\;\sigma^e_h\right)_{T^{-1}_{1/J_k}} + \sum_{e \in \mathcal{E}_{h,I}} \int_e C_{11}[\![v^h]\!] : [\![w^h]\!] \; d\Gamma\\
         =& \sum_{e \in \mathcal{E}_{h,I}}\sum_{\kappa \in \mathcal{T}_h} \left(r^e([\![w^h]\!]) + l^e\left(\vect{C_{12}}\cdot[\![w^h]\!]\right),\;\mathcal{G}^h\left(K\mathcal{G}^h(\nabla v^h) + Kr^e([\![v^h]\!]) + Kl^e\left(\vect{C_{12}}\cdot [\![v^h]\!]\right)\right)\right)_{T^{-1}_{1/J_k}}\\& + \sum_{e \in \mathcal{E}_{h,I}} \int_e C_{11}[\![v^h]\!] : [\![w^h]\!] \; d\Gamma.
    \end{split}
\end{equation}
Using \cref{eq:Gh} to the first term above yields
\begin{equation}
\label{eq:sigmahat_cdg}
    \begin{split}
        - \sum_{e \in \mathcal{E}_{h,I}} \int_e \hat{\sigma} : [\![w^h]\!] \; d\Gamma =& \sum_{\kappa \in \mathcal{T}_h} \int_\kappa \left(r_{\mathcal{I}}([\![w^h]\!]) + l_{\mathcal{I}}\left(\vect{C_{12}}\cdot[\![w^h]\!]\right)\right) : K\mathcal{G}^h(\nabla v^h) d\Omega \\
        &+ \sum_{e \in \mathcal{E}_{h,I}}\sum_{\kappa \in \mathcal{T}_h} \int_\kappa K \left(r^e([\![w^h]\!]) + l^e\left(\vect{C_{12}}\cdot[\![w^h]\!]\right)\right):\left( r^e([\![v^h]\!]) + l^e\left(\vect{C_{12}}\cdot [\![v^h]\!]\right)\right)d\Omega\\& + \sum_{e \in \mathcal{E}_{h,I}} \int_e C_{11}[\![v^h]\!] : [\![w^h]\!] \; d\Gamma.
    \end{split}
\end{equation}

Finally, substituting \cref{eq:sigmahat_cdg} into \cref{eq:BI_cdg_1}, one obtains
\begin{equation}
\label{eq:BI_cdg_final}
    \begin{split}
        B_{\mathcal{I}}^{v,\text{CDG-E}}(u^h,w^h) =& \sum_{\kappa \in \mathcal{T}_h} \int_\kappa K\mathcal{G}^h(\nabla w^h): \mathcal{G}^h(\nabla v^h) d\Omega + \sum_{\kappa \in \mathcal{T}_h} \int_\kappa  K\mathcal{G}^h(\nabla w^h):\left(r_{\mathcal{I}}([\![v^h]\!]) + l_{\mathcal{I}}(\vect{C_{12}}\cdot [\![v^h]\!])\right) d\Omega \\
    & + \sum_{\kappa \in \mathcal{T}_h} \int_\kappa K\mathcal{G}^h(\nabla v^h): \left(r_{\mathcal{I}}([\![w^h]\!]) + l_{\mathcal{I}}\left(\vect{C_{12}}\cdot[\![w^h]\!]\right)\right) d\Omega \\
        &+ \sum_{e \in \mathcal{E}_{h,I}}\sum_{\kappa \in \mathcal{T}_h} \int_\kappa K \left(r^e([\![w^h]\!]) + l^e\left(\vect{C_{12}}\cdot[\![w^h]\!]\right)\right):\left( r^e([\![v^h]\!]) + l^e\left(\vect{C_{12}}\cdot [\![v^h]\!]\right)\right)d\Omega\\& + \sum_{e \in \mathcal{E}_{h,I}} \int_e C_{11}[\![v^h]\!] : [\![w^h]\!] \; d\Gamma.
    \end{split}
\end{equation}

At the boundary, substituting \cref{eq:bc_cdg} into \cref{eq:B_gamma_final}, we get
\begin{equation}
\label{eq:B_gamma_cdg}
       B_{\Gamma}^{v,\text{CDG-E}}(u^h,w^h) =  \sum_{e \in \mathcal{E}_{h,\partial}}\sum_{\kappa \in \mathcal{T}_h}\int_\kappa K\mathcal{G}^h(\nabla w^h) : r^e\left((v_\Gamma(v^h) - v^h) \otimes n\right) \;d\Omega\;   -\sum_{e \in \mathcal{E}_{h,\partial}} \int_e \sigma_\Gamma(\mathcal{G}^h(K\mathcal{G}^h(\nabla v^h))) : [\![w^h]\!] \; d\Gamma,  
\end{equation}
which is the same as \cref{eq:B_gamma_br2}.

In order to prove entropy stability, we proceed as before and start by showing that $B^{v,\text{CDG-E}}_{\mathcal{I}}(u^h,v^h)\geq 0$.

\begin{theorem}
\label{thm:BI_cdg}
    For $C_{11}$ sufficiently large, we have $B^{v,\text{CDG-E}}_{\mathcal{I}}(u^h,v^h)\geq 0$.
\end{theorem}
\begin{proof}
    The proof is similar to that outlined in \cite{unified_viscous_dg} for the energy stability of Poisson's equation. We adapt the proof here for entropy stability and to obtain the value of $C_{11}$. Substituting $w^h=v^h$, we get
    \begin{equation}
\label{eq:BI_cdg_interimproof}
    \begin{split}
        B_{\mathcal{I}}^{v,\text{CDG-E}}(u^h,v^h) =& \sum_{\kappa \in \mathcal{T}_h} \int_\kappa K\mathcal{G}^h(\nabla v^h): \mathcal{G}^h(\nabla v^h) d\Omega + 2\sum_{\kappa \in \mathcal{T}_h} \int_\kappa  K\mathcal{G}^h(\nabla v^h):\left(r_{\mathcal{I}}([\![v^h]\!]) + l_{\mathcal{I}}(\vect{C_{12}}\cdot [\![v^h]\!])\right) d\Omega \\
        &+ \sum_{e \in \mathcal{E}_{h,I}}\sum_{\kappa \in \mathcal{T}_h} \int_\kappa K \left(r^e([\![v^h]\!]) + l^e\left(\vect{C_{12}}\cdot[\![v^h]\!]\right)\right):\left( r^e([\![v^h]\!]) + l^e\left(\vect{C_{12}}\cdot [\![v^h]\!]\right)\right)d\Omega\\& + \sum_{e \in \mathcal{E}_{h,I}} \int_e C_{11}[\![v^h]\!] : [\![v^h]\!] \; d\Gamma.
    \end{split}
\end{equation}
Since $K(v^h)$ is symmetric and positive semidefinite, $\sqrt{K}$ exists and is self-adjoint. For any $a,b \in \Sigma_{h,p}$, since $(\sqrt{K}a+\sqrt{K}b,\sqrt{K}a+\sqrt{K}b)_\kappa \geq 0$, we have $2(Ka,b)_{\kappa} \geq -(Ka,a)_{\kappa} - (Kb,b)_{\kappa} $. Using this inequality for the second term above results in:
    \begin{equation}
\label{eq:BI_cdg_interimproof2}
    \begin{split}
        B_{\mathcal{I}}^{v,\text{CDG-E}}(u^h,v^h) \geq& -\sum_{\kappa \in \mathcal{T}_h} \int_\kappa  K\left(r_{\mathcal{I}}([\![v^h]\!]) + l_{\mathcal{I}}(\vect{C_{12}}\cdot [\![v^h]\!])\right):\left(r_{\mathcal{I}}([\![v^h]\!]) + l_{\mathcal{I}}(\vect{C_{12}}\cdot [\![v^h]\!])\right) d\Omega \\
        &+ \sum_{e \in \mathcal{E}_{h,I}}\sum_{\kappa \in \mathcal{T}_h} \int_\kappa K \left(r^e([\![v^h]\!]) + l^e\left(\vect{C_{12}}\cdot[\![v^h]\!]\right)\right):\left( r^e([\![v^h]\!]) + l^e\left(\vect{C_{12}}\cdot [\![v^h]\!]\right)\right)d\Omega\\& + \sum_{e \in \mathcal{E}_{h,I}} \int_e C_{11}[\![v^h]\!] : [\![v^h]\!] \; d\Gamma.
    \end{split}
\end{equation}
If $N_f$ is the number of faces of an element, we have the following identity, adapted from \cite[Section 4.1]{unified_viscous_dg} using the fact that $r^e$ and $l^e$ are non-zero only on the elements that share the edge $e$:
\begin{equation}
\begin{split}
     \int_\kappa K \left(r_{\mathcal{I}}([\![v^h]\!]) + l_{\mathcal{I}}(\vect{C_{12}}\cdot [\![v^h]\!])\right)&:\left(r_{\mathcal{I}}([\![v^h]\!]) + l_{\mathcal{I}}(\vect{C_{12}}\cdot [\![v^h]\!])\right) d\Omega \leq \\
     & N_f\sum_{e\in \mathcal{E}_{h,I}}\int_\kappa  K\left(r^e([\![v^h]\!]) + l^e(\vect{C_{12}}\cdot [\![v^h]\!])\right):\left(r^e([\![v^h]\!]) + l^e(\vect{C_{12}}\cdot [\![v^h]\!])\right) d\Omega
     \end{split}
\end{equation}
    Therefore,
    \begin{equation}
\label{eq:BI_cdg_interimproof3}
    \begin{split}
        B_{\mathcal{I}}^{v,\text{CDG-E}}(u^h,v^h) \geq&\left(1-N_f\right) \sum_{e \in \mathcal{E}_{h,I}}\sum_{\kappa \in \mathcal{T}_h} \int_\kappa K \left(r^e([\![v^h]\!]) + l^e\left(\vect{C_{12}}\cdot[\![v^h]\!]\right)\right):\left( r^e([\![v^h]\!]) + l^e\left(\vect{C_{12}}\cdot [\![v^h]\!]\right)\right)d\Omega\\& + \sum_{e \in \mathcal{E}_{h,I}} \int_e C_{11}[\![v^h]\!] : [\![v^h]\!] \; d\Gamma \\
        \geq & \left(1-N_f\right) \|K\|_\infty \sum_{e \in \mathcal{E}_{h,I}}  \bigg\|r^e([\![v^h]\!]) + l^e\left(\vect{C_{12}}\cdot[\![v^h]\!]\right)\bigg\|_{L^2(\Omega)}^2 +   \sum_{e \in \mathcal{E}_{h,I}} C_{11} \bigg\|[\![v^h]\!]\bigg\|_{L^2(e)}^2,
    \end{split}
\end{equation}
where $\|K\|_\infty = \sup_{x \in \Omega} \{\|K(v^h(x)) \| \}$ and, for any valid mesh, $1-N_f<0$.
From \cite[Lemma 2]{brezzi_dg}, and the fact that the weight-adjusted norm is bounded from above and below by the usual $L^2-$ norm \cite{chan_weightadjusted_math}, we have the following inequalities for any edge function $\phi \in [L^2(e)]^{s\times d}$ and $q \in [L^2(e)]^s$:
\begin{equation}
\begin{aligned}
     &\|r^e(\phi)\|_{L^2(\Omega)}^2 \leq \frac{C_r}{h_e} \|\phi\|_{L^2(e)}^2, \\
      &\|l^e(q)\|_{L^2(\Omega)}^2 \leq \frac{C_l}{h_e} \|q\|_{L^2(e)}^2,
\end{aligned}
\end{equation}
where $h_e$ is the characteristic length of the edge $e$ and the constants $C_r$ and $C_l$ depend only on the number of faces per element $N_f$, the angle between the element faces and the polynomial degree \cite{brezzi_dg, unified_viscous_dg}. Therefore, we have
\begin{equation}
\label{eq:inequality_lifting}
    \bigg\|r^e([\![v^h]\!]) + l^e\left(\vect{C_{12}}\cdot[\![v^h]\!]\right)\bigg\|_{L^2(\Omega)}^2 \leq \left(\frac{C_r + \|\vect{C_{12}}\| C_l}{h_e} \right)\bigg\|[\![v^h]\!]\bigg\|_{L^2(e)}^2 = \frac{C^*}{h_e} \bigg\|[\![v^h]\!]\bigg\|_{L^2(e)}^2, 
\end{equation}
where $C^* = C_r + \|\vect{C_{12}}\| C_l$ depends on $\vect{C_{12}}$, $N_f$, the angle between the element faces and the polynomial degree.

Substituting \cref{eq:inequality_lifting} into \cref{eq:BI_cdg_interimproof3} yields
    \begin{equation}
\label{eq:BI_cdg_interimproof4}
    \begin{split}
        B_{\mathcal{I}}^{v,\text{CDG-E}}(u^h,v^h) \geq &  \sum_{e \in \mathcal{E}_{h,I}} \left(\frac{\left(1-N_f\right) \|K\|_\infty C^*}{h_e}  + C_{11}\right) \bigg\|[\![v^h]\!]\bigg\|_{L^2(e)}^2
    \end{split}
\end{equation}
Hence, if $C_{11} \geq \frac{\left(N_f-1\right) \|K\|_\infty C^*}{h_e}$, we have $B_{\mathcal{I}}^{v,\text{CDG-E}}(u^h,v^h) \geq 0$.

\end{proof}

\begin{theorem}
\label{thm:B_gamma_cdg}
   Denote $\tilde{\sigma}^h = \mathcal{G}^h(K\mathcal{G}^h(\nabla v^h))$. Let
    \begin{equation}
        v_\Gamma(v^h) = \frac{v^+(v^h)+v^h}{2}
    \end{equation}
    and 
    \begin{equation}
        \sigma_\Gamma(\tilde{\sigma}^h) = \frac{\sigma^+(\tilde{\sigma}^h) + \tilde{\sigma}^h}{2}.
    \end{equation}
    If the functions $v^+(\cdot)$ and $\sigma^+(\cdot)$ are chosen as specified in \cite{chan_entropystable_viscous}, then $\mathcal{B}^{v,\text{CDG-E}}_\Gamma(u^h,v^h)$ is a consistent discretization of the entropy boundary integral in \cref{eq:continuous_ns_entropy_inequality}, for all the boundary conditions presented in \cite{chan_entropystable_viscous}.
\end{theorem}
\begin{proof}
    The proof is the same as in \Cref{thm:B_gamma_br2}.
\end{proof}

\begin{theorem}
\label{thm:entropy_stability_cdg}
    Consider the CDG-E discretization scheme, with the penalization $C_{11}$ as given in the proof of \Cref{thm:BI_cdg}:
    \begin{equation}
    \sum_{\kappa \in \mathcal{T}_h}\left(\frac{\partial u^h}{\partial t},w^h\right)_{T^{-1}_{1/J_k}} + {B}^c(u^h,w^h) + B^{v,\text{CDG-E}}(u^h,w^h) = 0 \quad \forall w^h \in V_{h,p},
\end{equation}
where $B^{v,\text{CDG-E}}(u^h,w^h) = B_{\mathcal{I}}^{v,\text{CDG-E}}(u^h,w^h) + B_\Gamma^{v,\text{CDG-E}}(u^h,w^h)$, with $B_{\mathcal{I}}^{v,\text{CDG-E}}(u^h,w^h)$ given by \cref{eq:BI_cdg_final} and $B_\Gamma^{v,\text{CDG-E}}(u^h,w^h)$ given by \cref{eq:B_gamma_cdg}. If $\mathcal{B}^c(u^h,w^h)$ is an entropy-stable discretization of the convective term, the scheme is entropy-stable.
\end{theorem}
\begin{proof}
    Using \Cref{thm:BI_cdg} and \Cref{thm:B_gamma_cdg}, the proof is the same as in \Cref{thm:entropy_stability_br2} 
\end{proof}

\section{Numerical test cases}
We now present some numerical test cases to demonstrate the performance of both the schemes. For the convective residual $B^c(u^h,w^h)$, we use the Nonlinearly-Stable Flux Reconstruction (NSFR) scheme~\cite{alex_nsfr}. The entire scheme can be written as
\begin{equation}
\left(\vect{M}_m +\vect{K}_m\right)\frac{d}{d t} \hat{\vect{u}}_m(t)^T
      + \mathcal{B}^{c}(u^h,\chi) + \mathcal{B}^{v,\text{BR2-E/CDG-E}}(u^h,\chi)=\vect{0}^T,
  \end{equation}
  
  \noindent with the entropy-stable convective discretization \cite{alex_nsfr,chan_entropystable_dg}
  \begin{equation}
      \mathcal{B}^{c}(u^h,\chi) = \left[ \vect{\chi}(\vect{\xi}_v^r)^T\:  \vect{\chi}(\vect{\xi}_f^r)^T\right]
       \left[ \left(\Tilde{\vect{Q}}-\Tilde{\vect{Q}}^T\right) \odot \vect{F}_{m}^r \right]\vect{1}^T
      + \sum_{f=1}^{N_f}\sum_{k=1}^{N_{fp}}
     \vect{\chi}(\vect{\xi}_{fk}^r)^T W_{fk} \hat{\vect{n}}^r \cdot
     \vect{f}_m^{*,r}, 
  \end{equation}
  and
  \begin{equation}
      \Tilde{\vect{Q}} - \Tilde{\vect{Q}}^T = 
      \begin{bmatrix}
          \vect{W}\nabla^r\vect{\chi}(\vect{\xi}_v^r) - \nabla^r\vect{\chi}(\vect{\xi}_v^r) ^T \vect{W}
          &  \sum_{f=1}^{N_f} \vect{\chi}(\vect{\xi}_f^r)^T \vect{W}_f \text{diag}(\hat{\vect{n}}^r_f) \\
      -  \sum_{f=1}^{N_f}  \vect{W}_f \text{diag}(\hat{\vect{n}}^r_f) \vect{\chi}(\vect{\xi}_f^r) & \vect{0} 
      \end{bmatrix}.
  \end{equation}
In this work, we set the NSFR filter matrix to be $\vect{K}_m = 0$. The discrete form of the viscous term in an element can be written as (cf. \cref{eq:Bv4}):
\begin{equation}
\begin{aligned}
        \mathcal{B}^{v,\text{BR2-E/CDG-E}}(u^h,\chi) = \;& \vect{\mathcal{G}^h}\left(\nabla \vect{\chi}\right)(\vect{\xi}_{v}^r)^T \vect{J}(\vect{\xi}_{v}^r) \vect{W} \vect{K} \vect{\theta^h}(\vect{\xi}_{v}^r) - \sum_{f=1}^{N_f}\sum_{k=1}^{N_{fp}}
     \vect{\chi}(\vect{\xi}_{fk}^r)^T J_{fk} W_{fk} \hat{\vect{n}} \cdot
     \vect{\hat{\sigma}},
\end{aligned}
\end{equation}
where
\begin{equation}
    \vect{\theta^h}(\vect{\xi}_{v}^r) =  \vect{\mathcal{G}^h}\left(\nabla \vect{v^h}\right)(\vect{\xi}_{v}^r) + l_{\mathcal{I}}(\{\vect{v^h}\} - \hat{\vect{v}})(\vect{\xi}_{v}^r) + r_{\mathcal{I}}([\![\vect{v^h}]\!])(\vect{\xi}_{v}^r) - r_{\Gamma}((\vect{\hat{v} - v^h})\vect{\hat{n}})(\vect{\xi}_v^r), 
\end{equation}
with $r_{\mathcal{I}}$, $r_\Gamma$ and $l_{\mathcal{I}}$ as defined in \Cref{sec:weight_adjusted} and $\vect{\hat{v}}$ and $\vect{\hat{\sigma}}$ depending on the specific fluxes of BR2-E or CDG-E. 

The first test case is a three-dimensional lid-driven cavity that is used to check entropy stability, and the second case is a two-dimensional channel flow that is used to check the orders of convergence of the scheme. Test cases were run using our in-house code, Parallel High-Order Library for PDEs (PHiLiP) \cite{doug}, which uses deal.II \cite{dealii}.

\subsection{Three-dimensional lid-driven cavity}
The lid-driven cavity test case has been used in the literature to test the imposition of viscous boundary conditions. A two-dimensional version of the present test case was used in \cite{chan_entropystable_viscous} to test the entropy stability for the BR1 type scheme. We run the same case in three dimensions to test the entropy stability of BR2-E and CDG-E, with the derived modified boundary conditions.

\begin{figure}[H]

\begin{subfigure}{0.44\textwidth}
\centering
\includegraphics[width=1.0\linewidth]{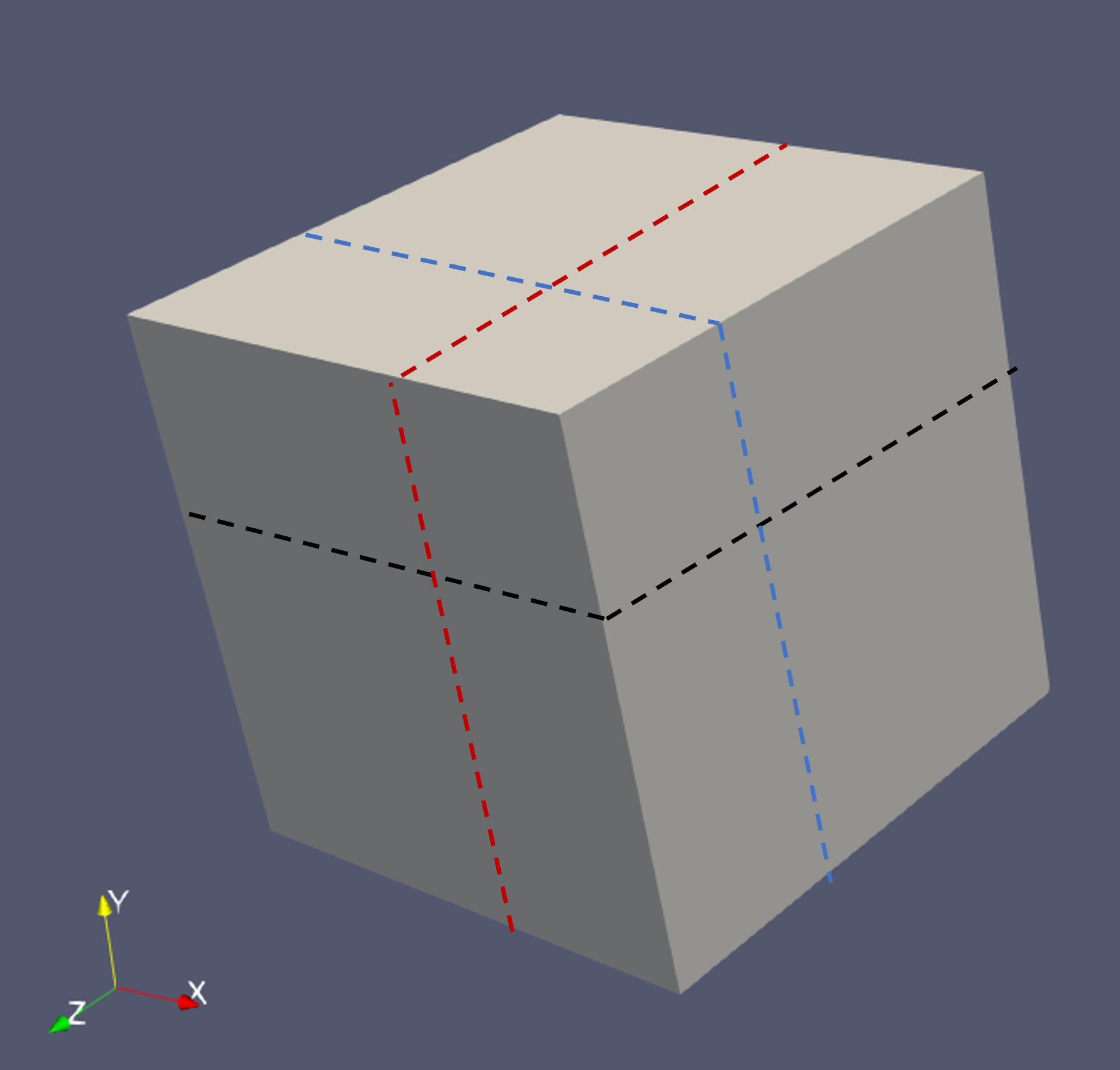} 
\caption{Cross-sections to visualize the grid in (b).}
\label{fig:subim1}
\end{subfigure}
\hfill
\begin{subfigure}{0.41\textwidth}
\centering
\includegraphics[width=1.0\linewidth]{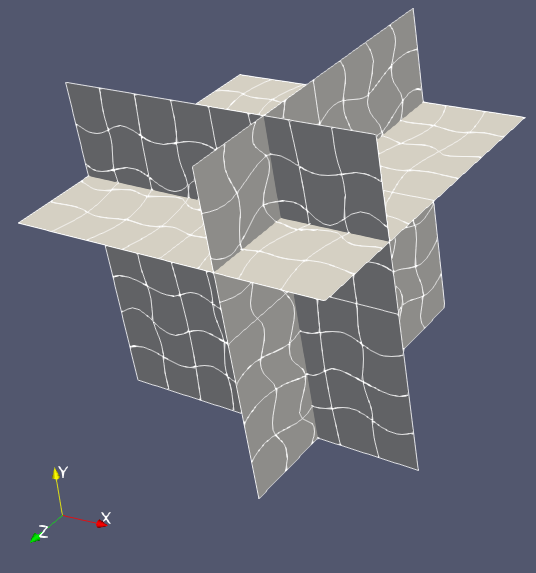}
\caption{Grid in various planes.}
\label{fig:subim2}
\end{subfigure}

\caption{Cross-sections of the cube domain $[-1,1]^3$ of the non-symmetric curvilinear grid \cite{alex_nsfr}.}
\label{fig:warpedgrid}
\end{figure}

The fluid is in a box $[-1\;1]^3$, with stationary walls at the sides and bottom, while the top wall moves at a speed of 1 unit parallel to the wall. The initial conditions are 
\begin{equation}
\begin{aligned}
    &\rho = 1,\\
    &u_1 = 0,\\
    &u_2 = 0,\\
    &u_3 = 0,\\
    &p = \frac{1}{M_\infty^2\gamma},
\end{aligned}
\end{equation}
where $\gamma = 1.4$, the reference mach number $M_\infty = 0.1$ and Reynold's number $Re_{\infty} = 1000$. The Navier-Stokes equations have been non-dimensionalized with respect to the free-stream conditions. We use no-slip boundary conditions at the bottom and side walls, while imposing the velocity of $[u_1, u_2, u_3] = [1,0,0]$ at the top wall. Both the BR2-E and the CDG-E viscous numerical fluxes were used to run the test case. The BR2-E factor $\eta_e$ was taken as $N_f+1$, and for CDG-E, we used $C_{11} = p(p+1)/{h_e}$, where $p$ is the polynomial degree. The CDG-E switch $\vect{C_{12}}$ on the face between cells $\kappa^+$ and $\kappa^-$ was chosen as $\frac{n^+}{2}$ if the cell number of $\kappa^+$ is less than that of $\kappa^-$. We used a polynomial degree of $p=3$ for both schemes. A non-symmetric curvilinear warped grid of degree $q=3$ was used, with the transformation \cite{alex_nsfr}:
\begin{equation}
    \begin{aligned}
        &\tilde{x} = x + \frac{1}{50} \sin\left(\pi x\right)\sin\left(\pi y\right)\sin\left(2\pi z\right), \\
         &\tilde{y} = y + \frac{1}{50} \sin\left(4\pi x\right)\sin\left(\pi y\right)\sin\left(3\pi z\right), \\
          &\tilde{z} = z + \frac{1}{50} \sin\left(2\pi x\right)\sin\left(5\pi y\right)\sin\left(\pi z\right),
    \end{aligned}
\end{equation}
where $(x,y,z) \in [-1,1]^3$ are the coordinates of a straight-edge grid and $(\tilde{x},\tilde{y},\tilde{z}) \in [-1,1]^3$ are the coordinates of the high-order transformed grid. This grid was used in the entropy-stability tests for the convective term in \cite{alex_nsfr}, and is obtained from a non-symmetric mapping that ensures that the metric terms do not cancel implicitly. \Cref{fig:warpedgrid} shows the cross-sections of the warped grid.

\begin{figure}
    \centering
    \includegraphics[trim={9cm 5cm 5cm 4.5cm},clip,width=0.4\linewidth]{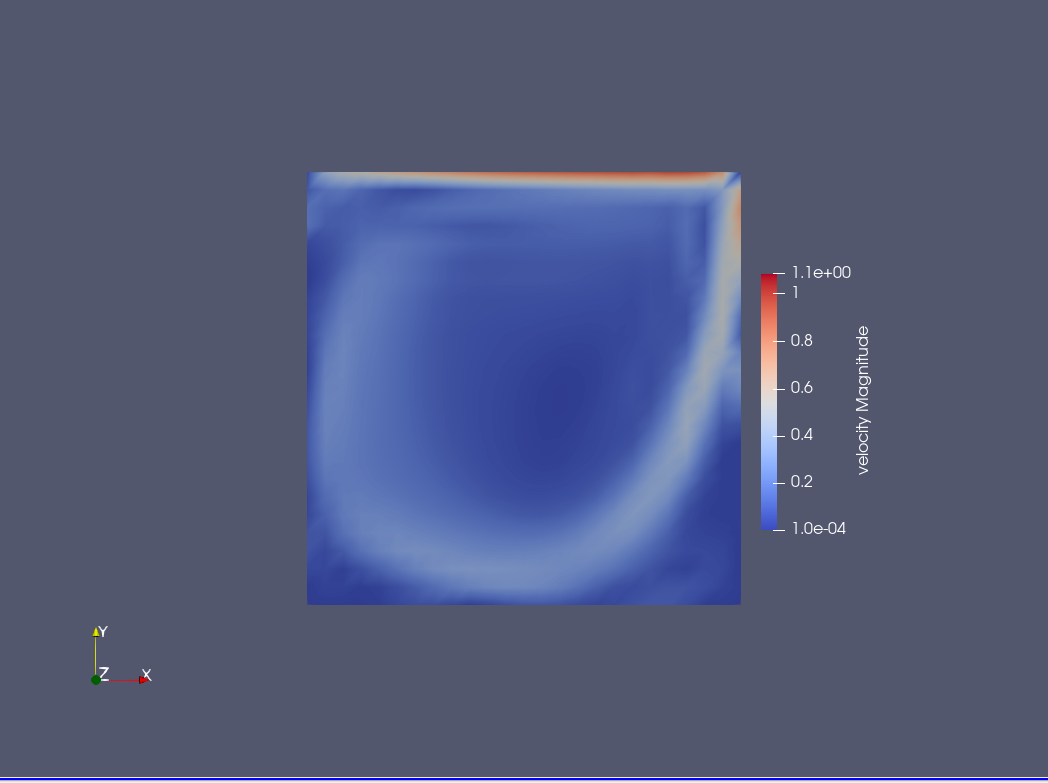}
    \caption{Velocity magnitude at $T=100$ with the BR2-E scheme.}
    \label{fig:br2_vel}
\end{figure}

\begin{figure}
    \centering
    \includegraphics[trim={9cm 5cm 5cm 4.5cm},clip,width=0.4\linewidth]{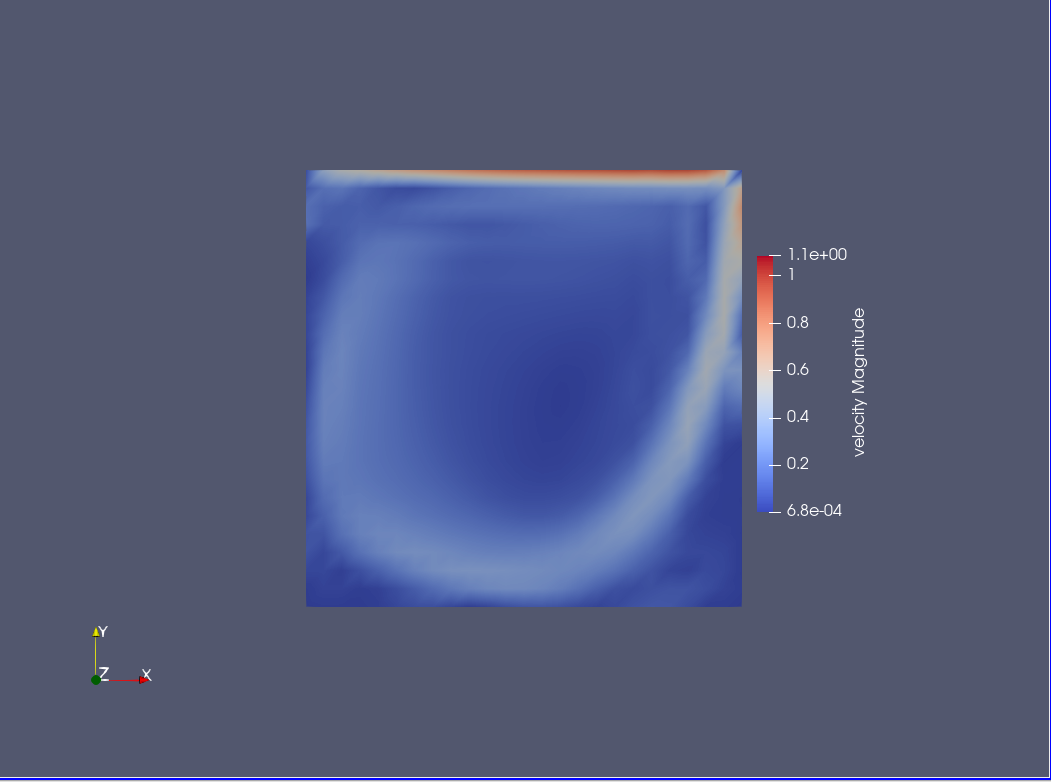}
    \caption{Velocity magnitude at $T=100$ with the CDG-E scheme.}
    \label{fig:cdg_vel}
\end{figure}

The convective terms were discretized using the NSFR scheme \cite{alex_nsfr}, which used the Chandrashekar two-point flux \cite{CH_flux} with Lax-Friedrichs dissipation. The Jacobian of the discretization was computed using automatic differentiation through \texttt{CodiPack} \cite{codipack}, and we used the third-order diagonally implicit Runge-Kutta \cite{dirk} with a constant time step of $0.5$. The simulation was run until a final time of $T = 100$. \Cref{fig:br2_vel} shows the final solution after $T=100$ using BR2-E, while \Cref{fig:cdg_vel} shows the final solution for CDG-E.

\begin{figure}[H]

\begin{subfigure}{1.0\textwidth}
\centering
\includegraphics[width=0.5\linewidth]{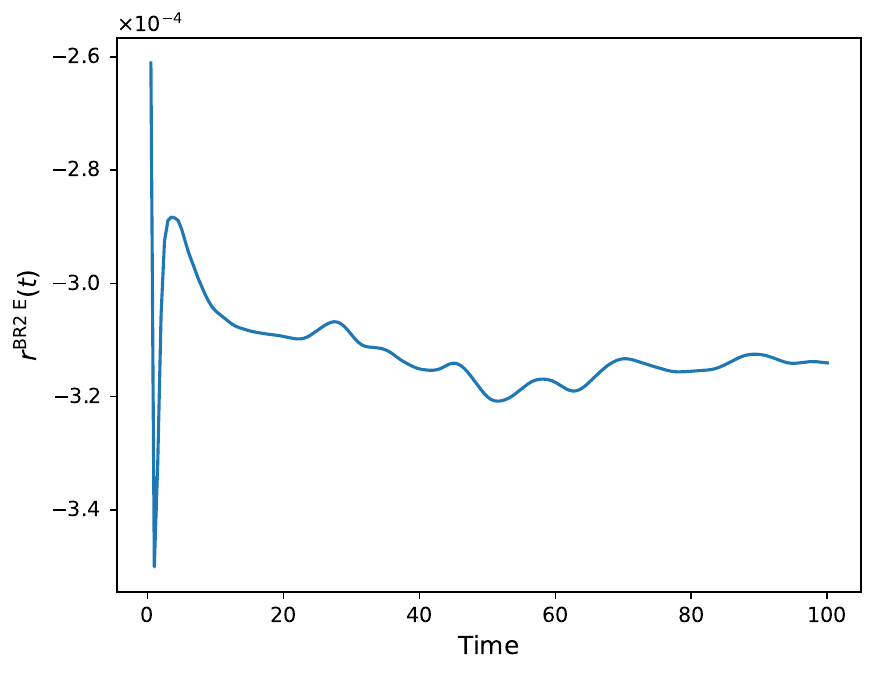} 
\caption{Plot of $r^{\text{BR2-E}}(t) = -\left(B^c(u^h,v^h) + B^{v,\text{BR2-E}}(u^h,v^h)\right) $ vs time.}
\label{fig:subim1}
\end{subfigure}
\begin{subfigure}{1.0\textwidth}
\centering
\includegraphics[width=0.5\linewidth]{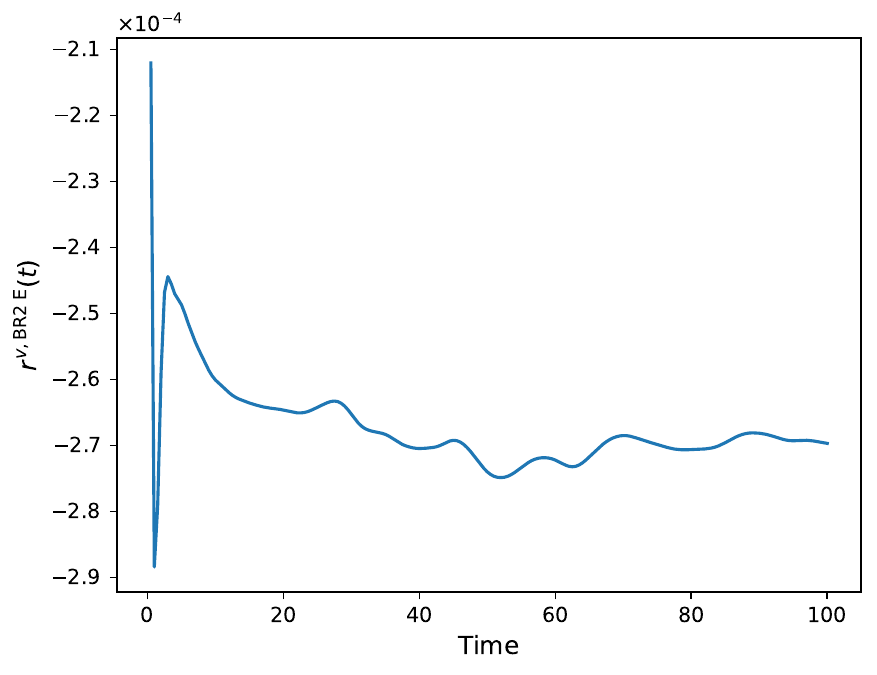}
\caption{Plot of $r^{v,\text{BR2-E}}(t) = - B^{v,\text{BR2-E}}(u^h,v^h) $ vs time. }
\label{fig:subim2}
\end{subfigure}

\caption{Plots of the inner product of the residuals with the entropy variable $v^h$ for the BR2-E discretization of the lid-driven cavity case.}
\label{fig:br2_v_times_rhs}
\end{figure}

\begin{figure}[H]

\begin{subfigure}{1.0\textwidth}
\centering
\includegraphics[width=0.5\linewidth]{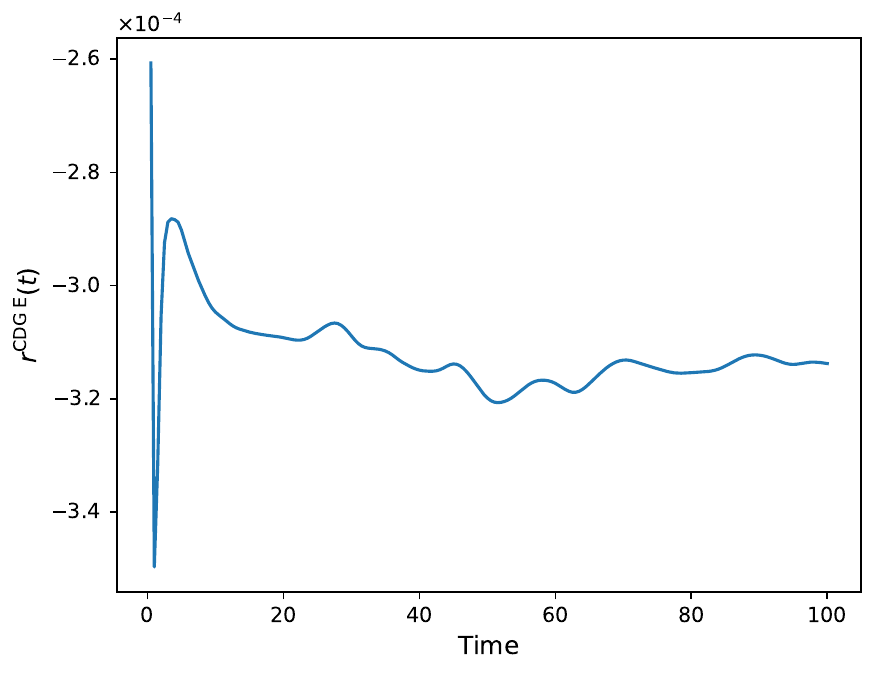} 
\caption{Plot of $r^{\text{CDG-E}}(t) = -\left(B^c(u^h,v^h) + B^{v,\text{CDG-E}}(u^h,v^h)\right) $ vs time.}
\label{fig:subim1}
\end{subfigure}
\begin{subfigure}{1.0\textwidth}
\centering
\includegraphics[width=0.5\linewidth]{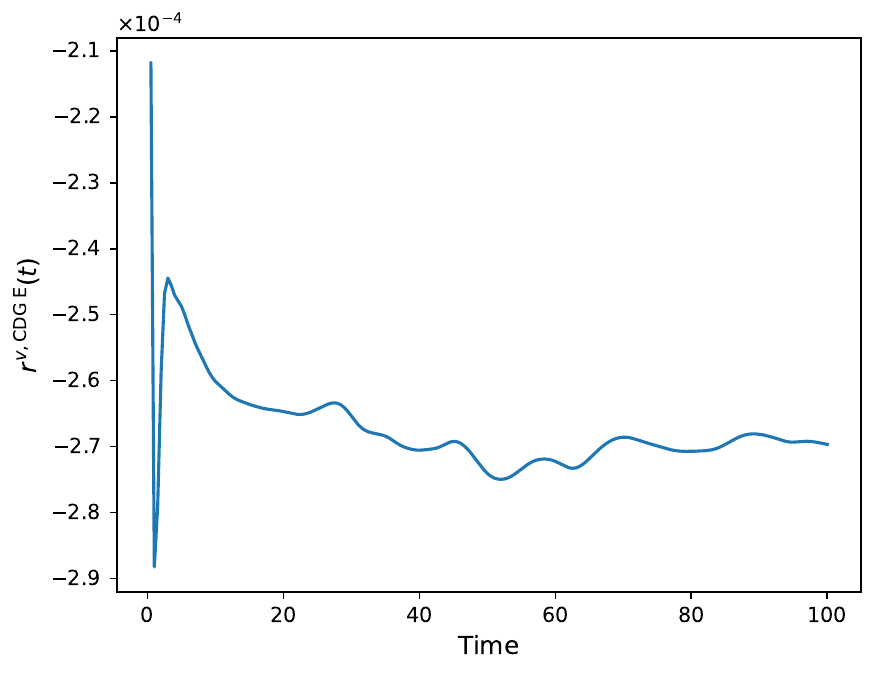}
\caption{Plot of $r^{v,\text{CDG-E}}(t) = - B^{v,\text{CDG-E}}(u^h,v^h) $ vs time. }
\label{fig:subim2}
\end{subfigure}

\caption{Plots of the inner product of the residuals with the entropy variable $v^h$ for the CDG-E discretization of the lid-driven cavity case.}
\label{fig:cdg_v_times_rhs}
\end{figure}

Next, we look at the entropy stability of the spatial discretization. To that end, we plot the inner products of $v^h$ with the right hand side, i.e.  $-\left(B^c(u^h,v^h) + B^{v,\text{BR2-E}}(u^h,v^h)\right)$ and $-B^{v,\text{BR2-E}}(u^h,v^h)$ in \Cref{fig:br2_v_times_rhs}. We also plot the inner products for CDG-E i.e $-\left(B^c(u^h,v^h) + B^{v,\text{CDG-E}}(u^h,v^h)\right)$ and $-B^{v,\text{CDG-E}}(u^h,v^h)$ in \Cref{fig:cdg_v_times_rhs}. Consistent with the proof, the spatial discretization of both the viscous terms as well as the entire residual is entropy-stable. We note that the differences in the inner-products of the BR2-E and CDG-E schemes are difficult to notice, indicating that most of the negative contribution is from the volume term.

\begin{figure}[H]

\begin{subfigure}{1.0\textwidth}
\centering
\includegraphics[width=0.5\linewidth]{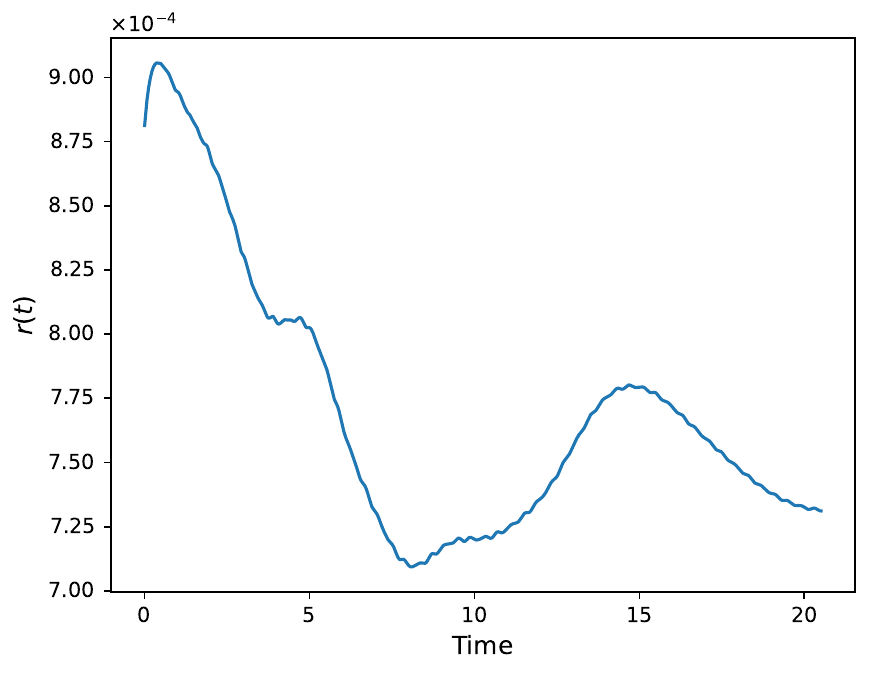} 
\caption{Plot of $r(t) = -\left(B^c(u^h,v^h) + B^{v}(u^h,v^h)\right) $ vs time.}
\label{fig:subim1}
\end{subfigure}
\begin{subfigure}{1.0\textwidth}
\centering
\includegraphics[width=0.5\linewidth]{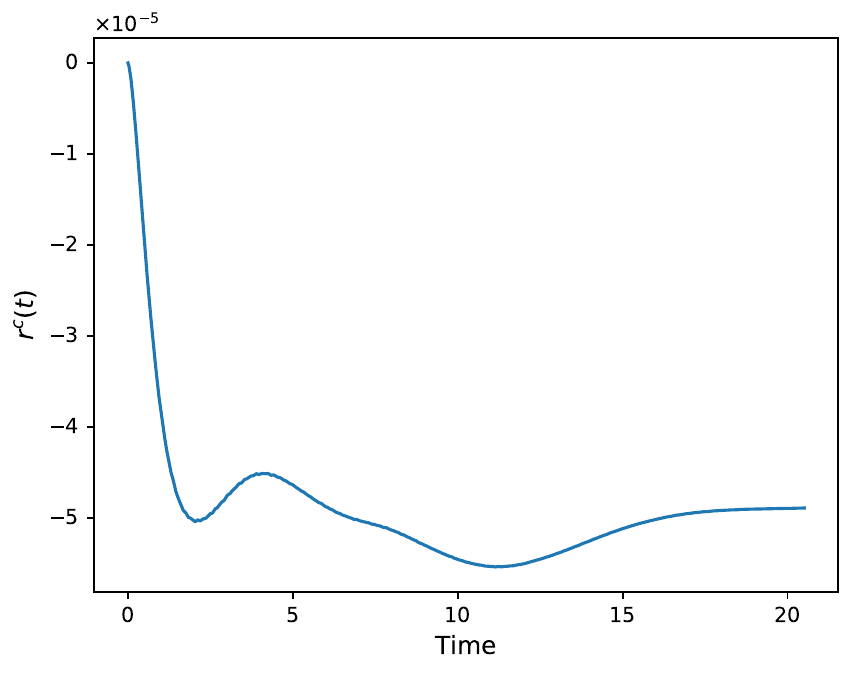}
\caption{Plot of $r^{c}(t) = - B^{c}(u^h,v^h) $ vs time. }
\label{fig:subim2}
\end{subfigure}

\caption{Plots of the inner product of the residuals with the entropy variable $v^h$. The convective term, $B^c(u^h,w^h)$, is discretized using the entropy-stable NSFR scheme \cite{alex_nsfr}, while the viscous term, $B^v(u^h,w^h)$, is discretized using the standard interior penalty method for DG \cite{julien_les,ip_arnold}. Note that $r(t)>0$, violating the entropy inequality.}
\label{fig:aux_v_times_rhs}
\end{figure}

Much focus in the literature has been on the entropy-stable discretization of the convective terms, with the assumption that since the viscous terms are dissipative by nature, they would likely not cause entropy violation. However, the present case of the lid-driven cavity shows the contrary. In fact, if we use entropy-stable discretization for just the convective term and use the standard DG discretization for the viscous term, we notice a violation of the entropy inequality. \Cref{fig:aux_v_times_rhs} plots the right-hand side of the entropy inequality, $r(t) = -\left(\mathcal{B}^c(u^h,v^h) + \mathcal{B}^v(u^h,v^h)\right)$, as well as the contribution to the right-hand side from the convective term $r^c(t) = -\mathcal{B}^c(u^h,v^h)$. As seen from the figure, the convective term, $r^c(t)$, is always negative as expected. However, the entropy violation from the dissipative term is significant enough for the net discretization to be entropy-violating.

\subsection{Periodic two-dimensional channel flow}

\begin{figure}[H]

\begin{subfigure}{1.0\textwidth}
\centering
\includegraphics[width=0.5\linewidth]{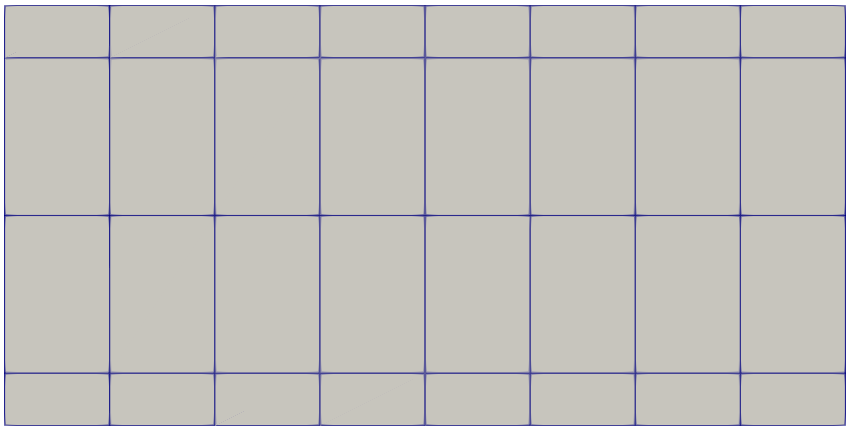} 
\caption{Coarsest grid for the 2D channel flow test case.}
\label{fig:subim1}
\end{subfigure}
\begin{subfigure}{1.0\textwidth}
\centering
\includegraphics[width=0.5\linewidth]{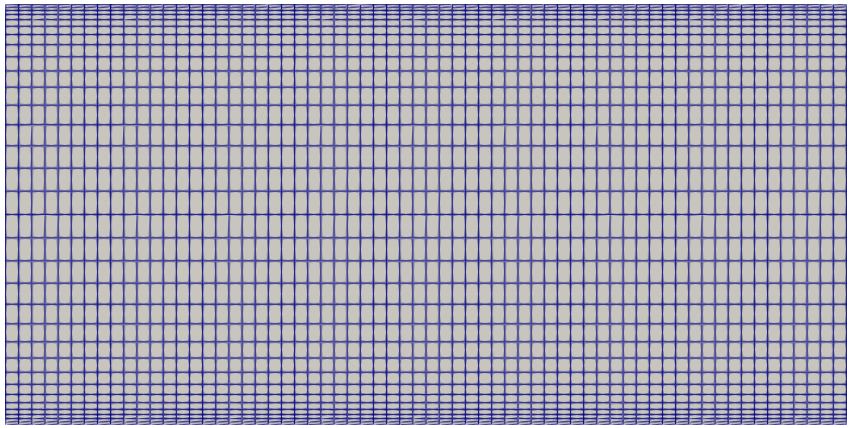}
\caption{Finest grid used for the 2D channel flow test case.}
\label{fig:subim2}
\end{subfigure}

\caption{Coarsest and finest grids used in the sequence for the channel flow test case.}
\label{fig:coarse_fine_channelflow}
\end{figure}

We run the viscous discretization schemes on a test case of the two-dimensional channel flow. The domain is rectangular $[-2,2] \times [-1,1]$, with solid walls at the top and bottom and periodicity in the $x$ direction. No-slip boundary conditions are imposed at the walls. We run the cases until a final time of $T=0.5$. This test case is the same as in \cite{chan_entropystable_viscous}, with the initial conditions:
\begin{equation}
    \begin{aligned}
        &\rho = 1,\\
        &u_1 = \frac{1}{10} \sin\left(\frac{\pi x}{2}\right)\cos\left(\frac{\pi y}{2}\right), \\
        &u_2 = \frac{1}{10} \cos\left(\frac{\pi x}{2}\right)\sin\left(\pi y\right),\\
        & p = \frac{1}{M_\infty^2 \gamma},
    \end{aligned}
\end{equation}
where $M_\infty = 0.3$ and $Re_\infty = 50$. We use $\eta_e = N_f + 1$ for BR2-E and $C_{11} = 5p(p+1)/h_e$ for CDG-E. We run both schemes for the polynomial degrees $p=1,2,3$ and $4$, and four grids with 32, 128, 512 and 2048 elements. Similar to \cite{chan_entropystable_viscous}, we transform the $y$-coordinates of a uniform mesh using the transformation $\tilde{y} = y + 0.25\sin(\pi y)$ to obtain a mesh that's fine near the wall. \Cref{fig:coarse_fine_channelflow} shows the coarsest and the finest mesh in the sequence.

We are interested in the convergence of the error in velocity at the wall, which is expected to be zero. To that end, we compute the $L_2$ error $e = \sqrt{ \int_{\partial \Omega_\text{wall}} u_1^2 + u_2^2 d\Gamma}$ for each grid and polynomial degree. 

The convective terms are discretized using the NSFR scheme \cite{alex_nsfr}, and we use third-order diagonally implicit Runge-Kutta \cite{dirk} with a fixed time step of $\Delta t = 0.01$ across all the grids and polynomial degrees. The simulations were run until a final time of $T=0.5$ before computing the error. \Cref{table:br2_convergence} shows the convergence of the error for the entropy-stable BR2-E scheme, while \Cref{table:cdg_convergence} shows the convergence of the entropy-stable CDG-E scheme. Both the schemes yield the expected $\mathcal{O}\left(h^{p+1}\right)$ rate of convergence.

\begin{table}[pos=h]
    \centering
    \caption{Convergence of the error for the BR2-E scheme}
    \label{table:br2_convergence}
    \begin{tabular}{|c|c|c|c|c|}
        \hline
        \textbf{Polynomial degree} & \textbf{Number of cells} & \textbf{Degrees of freedom} & \textbf{$L_2$ error} & \textbf{Rate of convergence}\\
        \hline
        1 & 32 & 512 & $1.46\times 10^{-1}$& -\\
        \hline
        1 & 128 & 2048 & $8.27\times 10^{-2}$& 0.8\\
        \hline
        1 & 512 & 8192  & $1.91\times10^{-2}$&2.1\\
        \hline
        1 & 2048 & 32768 & $3.88\times10^{-3}$&2.3\\
        \hline
        \hline
        2 & 32 & 1152 &$1.09\times10^{-1}$ & -\\
        \hline
        2 & 128 & 4608 &$8.35\times10^{-4}$ &7.0\\
        \hline
        2 & 512 & 18432 &$1.08\times 10^{-4}$ &2.9\\
        \hline
        2 & 2048 & 73728 &$1.70\times10^{-5}$ &2.7\\
         \hline
        \hline
        3 & 32 & 2048 &$3.54\times10^{-3}$ & -\\
        \hline
        3 & 128 & 8192 &$3.39\times10^{-4}$ &3.4\\
        \hline
        3 & 512 & 32768 &$1.59\times10^{-5}$ &4.4\\
        \hline
        3 & 2048 & 131072 &$9.10\times10^{-7}$ &4.1\\
         \hline
        \hline
        4 & 32 & 3200 &$3.60\times10^{-3}$ &-\\
        \hline
        4 & 128 & 12800 & $3.81\times10^{-5}$&6.6\\
        \hline
        4 & 512 & 51200 &$4.97\times10^{-7}$ &6.3\\
        \hline
        4 & 2048 & 204800 &$1.28\times10^{-8}$ &5.3\\

        \hline
    \end{tabular}
\end{table}

\begin{table}[pos=h]
    \centering
    \caption{Convergence of the error for the CDG-E scheme}
    \label{table:cdg_convergence}
    \begin{tabular}{|c|c|c|c|c|}
        \hline
        \textbf{Polynomial degree} & \textbf{Number of cells} & \textbf{Degrees of freedom} & \textbf{$L_2$ error} & \textbf{Rate of convergence}\\
        \hline
        1 & 32 & 512 & $1.47\times 10^{-1}$& -\\
        \hline
        1 & 128 & 2048 & $8.31\times 10^{-2}$& 0.8\\
        \hline
        1 & 512 & 8192  & $1.88\times10^{-2}$&2.1\\
        \hline
        1 & 2048 & 32768 & $3.86\times10^{-3}$&2.2\\
        \hline
        \hline
        2 & 32 & 1152 &$1.10\times10^{-1}$ & -\\
        \hline
        2 & 128 & 4608 &$1.09\times10^{-3}$ &6.7\\
        \hline
        2 & 512 & 18432 &$1.12\times 10^{-4}$ &3.3\\
        \hline
        2 & 2048 & 73728 &$1.90\times10^{-5}$ &2.7\\
         \hline
        \hline
        3 & 32 & 2048 &$4.04\times10^{-3}$ & -\\
        \hline
        3 & 128 & 8192 &$3.75\times10^{-4}$ &3.4\\
        \hline
        3 & 512 & 32768 &$1.59\times10^{-5}$ &4.6\\
        \hline
        3 & 2048 & 131072 &$9.07\times10^{-7}$ &4.1\\
         \hline
        \hline
        4 & 32 & 3200 &$3.63\times10^{-3}$ &-\\
        \hline
        4 & 128 & 12800 & $3.90\times10^{-5}$&6.5\\
        \hline
        4 & 512 & 51200 &$5.87\times10^{-7}$ &6.0\\
        \hline
        4 & 2048 & 204800 &$1.36\times10^{-8}$ &5.4\\

        \hline
    \end{tabular}
\end{table}

\begin{figure}[H]

\begin{subfigure}{1.0\textwidth}
\centering
\includegraphics[width=0.5\linewidth]{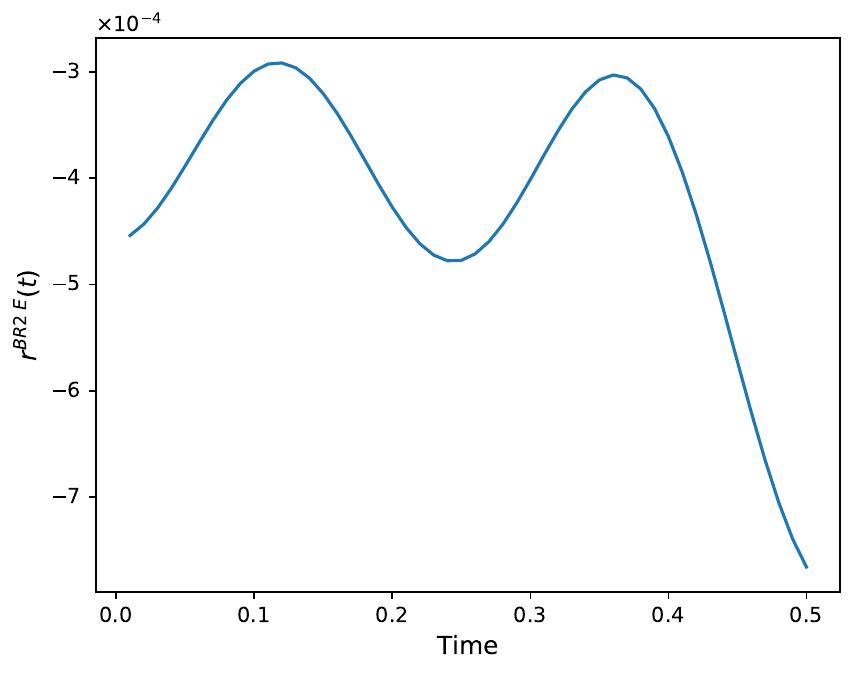} 
\caption{Plot of $r^{\text{BR2-E}}(t) = -\left(B^c(u^h,v^h) + B^{v,\text{BR2-E}}(u^h,v^h)\right) $ vs time.}
\label{fig:subim1}
\end{subfigure}
\begin{subfigure}{1.0\textwidth}
\centering
\includegraphics[width=0.5\linewidth]{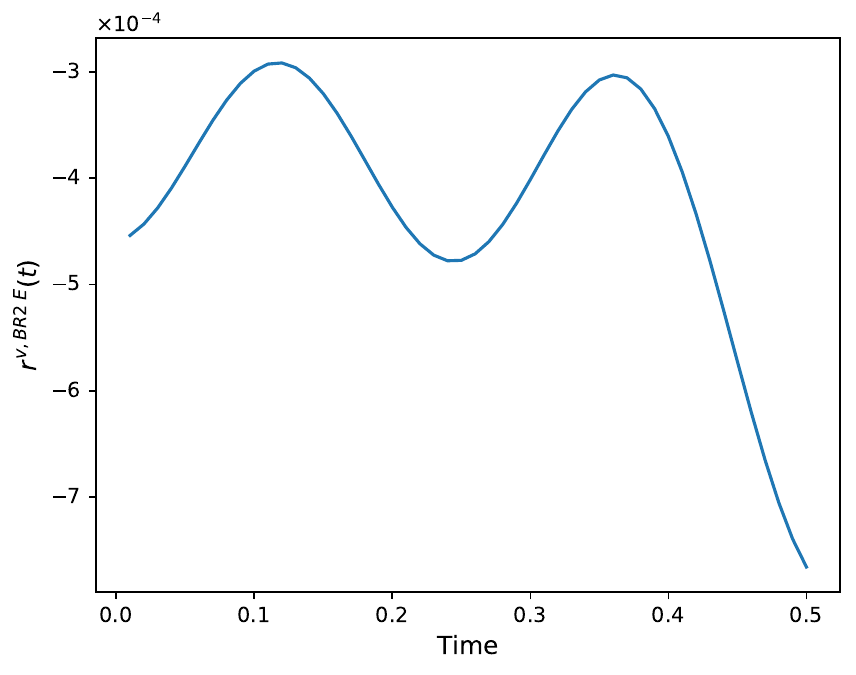}
\caption{Plot of $r^{v,\text{BR2-E}}(t) = - B^{v,\text{BR2-E}}(u^h,v^h) $ vs time. }
\label{fig:subim2}
\end{subfigure}

\caption{Plots of the inner product of the residuals with the entropy variable $v^h$ for the BR2-E discretization of the channel flow case.}
\label{fig:br2_v_times_rhs_channel}
\end{figure}

\begin{figure}[H]

\begin{subfigure}{1.0\textwidth}
\centering
\includegraphics[width=0.5\linewidth]{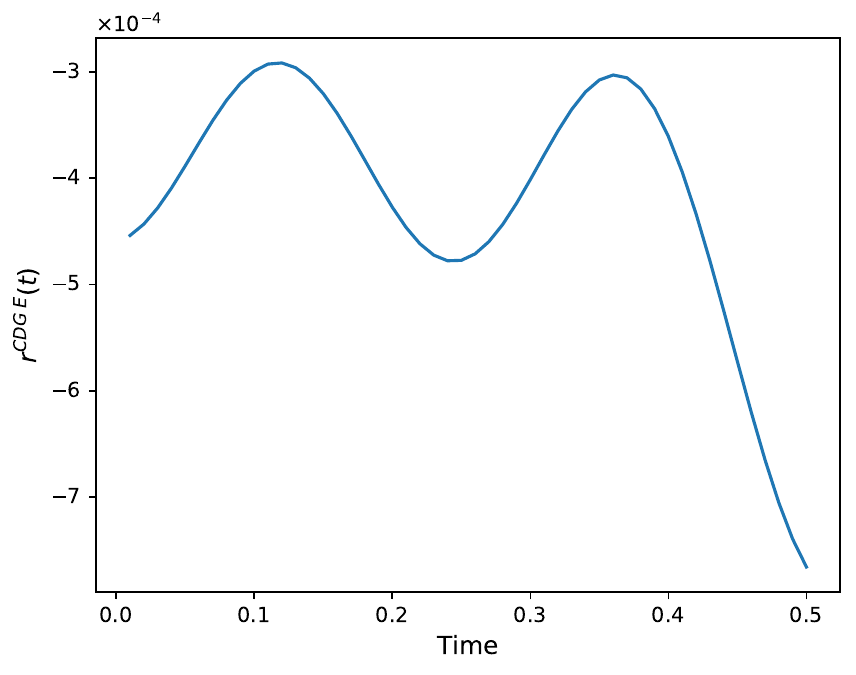} 
\caption{Plot of $r^{\text{CDG-E}}(t) = -\left(B^c(u^h,v^h) + B^{v,\text{CDG-E}}(u^h,v^h)\right) $ vs time.}
\label{fig:subim1}
\end{subfigure}
\begin{subfigure}{1.0\textwidth}
\centering
\includegraphics[width=0.5\linewidth]{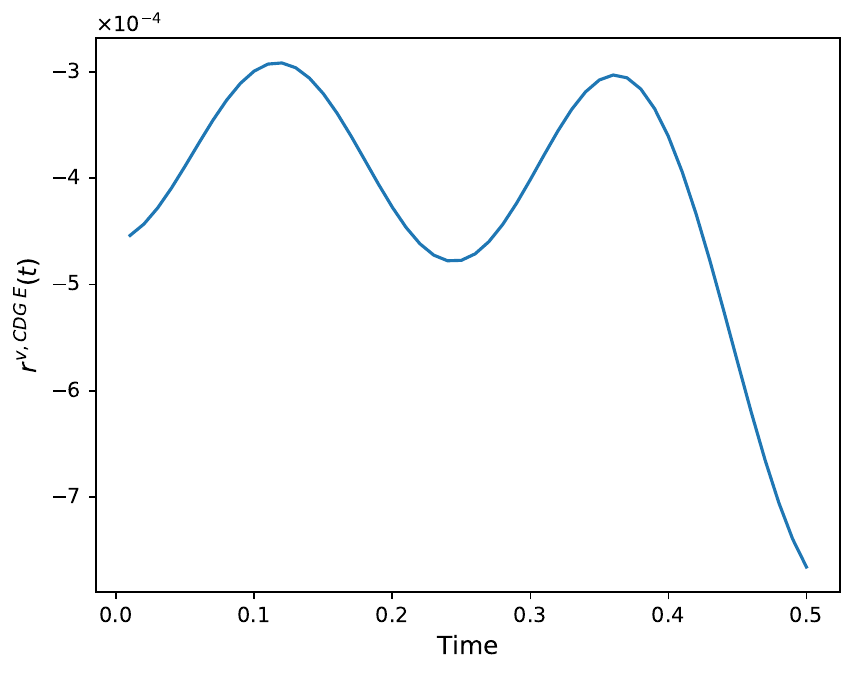}
\caption{Plot of $r^{v,\text{CDG-E}}(t) = - B^{v,\text{CDG-E}}(u^h,v^h) $ vs time. }
\label{fig:subim2}
\end{subfigure}

\caption{Plots of the inner product of the residuals with the entropy variable $v^h$ for the CDG-E discretization of the channel flow case.}
\label{fig:cdg_v_times_rhs_channel}
\end{figure}

 We also plot the evolution of the entropy residual for the finest grid with $p=4$. \Cref{fig:br2_v_times_rhs_channel} shows the entropy residual for BR2-E, while \Cref{fig:cdg_v_times_rhs_channel} shows the same for the CDG-E scheme. Consistent with the entropy-stability proofs, both the schemes have negative entropy residuals. 

\section{Conclusion}
The current work developed weight-adjusted formulations for the BR2 and CDG schemes, referred to as BR2-E and CDG-E, that are entropy-stable. The weight-adjusted formulation was used in the lifting and projection operators, as well as in the numerical fluxes, and entropy stability was proven for the two schemes. Another novelty of the current approach was the modification made to the numerical fluxes at the boundary that led to entropy stability for the boundary terms. We ran two- and three-dimensional test cases, demonstrating high-order convergence and entropy stability for both schemes. The primary motivation for the current work was the fact that the BR2 and CDG schemes have a compact stencil, as the residual on an element depends only on the solution in the element and in the immediate-adjacent elements. A compact stencil is useful because it allows efficient parallelizability and computation of the residual, while also permitting an easier formation of the Jacobian through automatic differentiation.

We note that the current work was on deriving and proving properties of the entropy-stable weight-adjusted formulations, along with the imposition of entropy-stable boundary fluxes. However, we leave the comparison between the performance of the various schemes for the future. To that end, it would be interesting to test the schemes on various turbulent flow cases, while noting the resolution of the energy spectrum and grid convergence of important outputs such as the lift and drag coefficients.

\printcredits

\bibliographystyle{model1-num-names}

\bibliography{cas-refs}

\begin{thebibliography}{28}
\expandafter\ifx\csname natexlab\endcsname\relax\def\natexlab#1{#1}\fi
\providecommand{\url}[1]{\texttt{#1}}
\providecommand{\href}[2]{#2}
\providecommand{\path}[1]{#1}
\providecommand{\DOIprefix}{doi:}
\providecommand{\ArXivprefix}{arXiv:}
\providecommand{\URLprefix}{URL: }
\providecommand{\Pubmedprefix}{pmid:}
\providecommand{\doi}[1]{\href{http://dx.doi.org/#1}{\path{#1}}}
\providecommand{\Pubmed}[1]{\href{pmid:#1}{\path{#1}}}
\providecommand{\bibinfo}[2]{#2}
\ifx\xfnm\relax \def\xfnm[#1]{\unskip,\space#1}\fi
\bibitem[{Hesthaven and Warburton(2008)}]{hesthaven2008nodal}
\bibinfo{author}{J.~S. Hesthaven}, \bibinfo{author}{T.~Warburton}, \bibinfo{title}{Nodal discontinuous Galerkin methods: Algorithms, Analysis, and Applications}, \bibinfo{publisher}{Springer}, \bibinfo{year}{2008}.
\bibitem[{Persson and Peraire(2006)}]{subcell_shock_capturing}
\bibinfo{author}{P.-O. Persson}, \bibinfo{author}{J.~Peraire},
\newblock \bibinfo{title}{Sub-cell shock capturing for discontinuous {Galerkin} methods},
\newblock in: \bibinfo{booktitle}{44th AIAA Aerospace Sciences Meeting and Exhibit}, \bibinfo{year}{2006}, p. \bibinfo{pages}{112}.
\bibitem[{Nordstr{\"o}m and Winters(2021)}]{filtering}
\bibinfo{author}{J.~Nordstr{\"o}m}, \bibinfo{author}{A.~R. Winters},
\newblock \bibinfo{title}{Stable filtering procedures for nodal discontinuous {G}alerkin methods},
\newblock \bibinfo{journal}{Journal of Scientific Computing} \bibinfo{volume}{87} (\bibinfo{year}{2021}) \bibinfo{pages}{17}.
\bibitem[{Cicchino and Nadarajah(2025)}]{alex_nsfr}
\bibinfo{author}{A.~Cicchino}, \bibinfo{author}{S.~Nadarajah},
\newblock \bibinfo{title}{Discretely nonlinearly stable weight-adjusted flux reconstruction high-order method for compressible flows on curvilinear grids},
\newblock \bibinfo{journal}{Journal of Computational Physics} \bibinfo{volume}{521} (\bibinfo{year}{2025}) \bibinfo{pages}{113532}.
\bibitem[{Chan(2018)}]{chan_entropystable_dg}
\bibinfo{author}{J.~Chan},
\newblock \bibinfo{title}{{On discretely entropy conservative and entropy stable discontinuous Galerkin methods}},
\newblock \bibinfo{journal}{Journal of Computational Physics} \bibinfo{volume}{362} (\bibinfo{year}{2018}) \bibinfo{pages}{346--374}.
\bibitem[{Lax(1971)}]{lax_entropy}
\bibinfo{author}{P.~Lax},
\newblock \bibinfo{title}{{Shock Waves and Entropy}},
\newblock in: \bibinfo{editor}{E.~H. Zarantonello} (Ed.), \bibinfo{booktitle}{Contributions to Nonlinear Functional Analysis}, \bibinfo{publisher}{Academic Press}, \bibinfo{year}{1971}, pp. \bibinfo{pages}{603--634}. \DOIprefix\doi{https://doi.org/10.1016/B978-0-12-775850-3.50018-2}.
\bibitem[{Brillon and Nadarajah(2026)}]{julien_les}
\bibinfo{author}{J.~Brillon}, \bibinfo{author}{S.~Nadarajah},
\newblock \bibinfo{title}{{Large eddy simulation using nonlinearly stable flux reconstruction}},
\newblock \bibinfo{journal}{Computers \& Fluids} \bibinfo{volume}{319} (\bibinfo{year}{2026}) \bibinfo{pages}{107232}.
\bibitem[{Zakerzadeh and May(2017)}]{May_viscous_entropystableDG}
\bibinfo{author}{M.~Zakerzadeh}, \bibinfo{author}{G.~May},
\newblock \bibinfo{title}{{Entropy stable discontinuous Galerkin scheme for the compressible Navier-Stokes equations}},
\newblock in: \bibinfo{booktitle}{55th AIAA Aerospace Sciences Meeting}, \bibinfo{year}{2017}, p. \bibinfo{pages}{0084}.
\bibitem[{Zhao et~al.(2027)Zhao, Zhao, Liu, and Hao}]{viscous_sipg_br1}
\bibinfo{author}{Z.~Zhao}, \bibinfo{author}{M.~Zhao}, \bibinfo{author}{Z.~Liu}, \bibinfo{author}{S.~Hao},
\newblock \bibinfo{title}{{Entropy stable DGSEM discretization of viscous terms in Navier-Stokes equations}},
\newblock \bibinfo{journal}{Journal of Computational Physics} \bibinfo{volume}{569} (\bibinfo{year}{2027}) \bibinfo{pages}{115430}.
\bibitem[{Chan et~al.(2022)Chan, Lin, and Warburton}]{chan_entropystable_viscous}
\bibinfo{author}{J.~Chan}, \bibinfo{author}{Y.~Lin}, \bibinfo{author}{T.~Warburton},
\newblock \bibinfo{title}{{Entropy stable modal discontinuous Galerkin schemes and wall boundary conditions for the compressible Navier-Stokes equations}},
\newblock \bibinfo{journal}{Journal of Computational Physics} \bibinfo{volume}{448} (\bibinfo{year}{2022}) \bibinfo{pages}{110723}.
\bibitem[{Sv{\"a}rd and {\"O}zcan(2014)}]{svard_entropystablebc1}
\bibinfo{author}{M.~Sv{\"a}rd}, \bibinfo{author}{H.~{\"O}zcan},
\newblock \bibinfo{title}{{Entropy-stable schemes for the Euler equations with far-field and wall boundary conditions}},
\newblock \bibinfo{journal}{Journal of Scientific Computing} \bibinfo{volume}{58} (\bibinfo{year}{2014}) \bibinfo{pages}{61--89}.
\bibitem[{Svärd(2021)}]{svard_entropystablebc2}
\bibinfo{author}{M.~Svärd},
\newblock \bibinfo{title}{{Entropy stable boundary conditions for the Euler equations}},
\newblock \bibinfo{journal}{Journal of Computational Physics} \bibinfo{volume}{426} (\bibinfo{year}{2021}) \bibinfo{pages}{109947}.
\bibitem[{Chan and Wilcox(2019)}]{chan_weightadjusted_jcp}
\bibinfo{author}{J.~Chan}, \bibinfo{author}{L.~C. Wilcox},
\newblock \bibinfo{title}{{On discretely entropy stable weight-adjusted discontinuous Galerkin methods: curvilinear meshes}},
\newblock \bibinfo{journal}{Journal of Computational Physics} \bibinfo{volume}{378} (\bibinfo{year}{2019}) \bibinfo{pages}{366--393}.
\bibitem[{Masatsuka(2009)}]{masatsuka}
\bibinfo{author}{K.~Masatsuka}, \bibinfo{title}{{I do like CFD, VOL.1, Second Edition}}, \bibinfo{year}{2009}.
\bibitem[{Hughes et~al.(1986)Hughes, Franca, and Mallet}]{K_matrix_3D}
\bibinfo{author}{T.~Hughes}, \bibinfo{author}{L.~Franca}, \bibinfo{author}{M.~Mallet},
\newblock \bibinfo{title}{{A new finite element formulation for computational fluid dynamics: I. Symmetric forms of the compressible Euler and Navier-Stokes equations and the second law of thermodynamics}},
\newblock \bibinfo{journal}{Computer Methods in Applied Mechanics and Engineering} \bibinfo{volume}{54} (\bibinfo{year}{1986}) \bibinfo{pages}{223--234}.
\bibitem[{Mock(1980)}]{mock}
\bibinfo{author}{M.~S. Mock},
\newblock \bibinfo{title}{Systems of conservation laws of mixed type},
\newblock \bibinfo{journal}{Journal of Differential equations} \bibinfo{volume}{37} (\bibinfo{year}{1980}) \bibinfo{pages}{70--88}.
\bibitem[{Dalcin et~al.(2019)Dalcin, Rojas, Zampini, {Del Rey Fernández}, Carpenter, and Parsani}]{dalcin}
\bibinfo{author}{L.~Dalcin}, \bibinfo{author}{D.~Rojas}, \bibinfo{author}{S.~Zampini}, \bibinfo{author}{D.~C. {Del Rey Fernández}}, \bibinfo{author}{M.~H. Carpenter}, \bibinfo{author}{M.~Parsani},
\newblock \bibinfo{title}{{Conservative and entropy stable solid wall boundary conditions for the compressible Navier–Stokes equations: Adiabatic wall and heat entropy transfer}},
\newblock \bibinfo{journal}{Journal of Computational Physics} \bibinfo{volume}{397} (\bibinfo{year}{2019}) \bibinfo{pages}{108775}.
\bibitem[{Chan et~al.(2017)Chan, Hewett, and Warburton}]{chan_weightadjusted_math}
\bibinfo{author}{J.~Chan}, \bibinfo{author}{R.~J. Hewett}, \bibinfo{author}{T.~Warburton},
\newblock \bibinfo{title}{{Weight-Adjusted Discontinuous Galerkin Methods: Wave Propagation in Heterogeneous Media}},
\newblock \bibinfo{journal}{SIAM Journal on Scientific Computing} \bibinfo{volume}{39} (\bibinfo{year}{2017}) \bibinfo{pages}{A2935--A2961}.
\bibitem[{Bassi and Rebay(2000)}]{bassi_rebay2}
\bibinfo{author}{F.~Bassi}, \bibinfo{author}{S.~Rebay},
\newblock \bibinfo{title}{{GMRES discontinuous Galerkin solution of the compressible Navier-Stokes equations}},
\newblock in: \bibinfo{booktitle}{Discontinuous Galerkin methods: Theory, Computation and Applications}, \bibinfo{publisher}{Springer}, \bibinfo{year}{2000}, pp. \bibinfo{pages}{197--208}.
\bibitem[{Arnold et~al.(2002)Arnold, Brezzi, Cockburn, and Marini}]{unified_viscous_dg}
\bibinfo{author}{D.~N. Arnold}, \bibinfo{author}{F.~Brezzi}, \bibinfo{author}{B.~Cockburn}, \bibinfo{author}{L.~D. Marini},
\newblock \bibinfo{title}{{Unified analysis of discontinuous Galerkin methods for elliptic problems}},
\newblock \bibinfo{journal}{SIAM Journal on Numerical Analysis} \bibinfo{volume}{39} (\bibinfo{year}{2002}) \bibinfo{pages}{1749--1779}.
\bibitem[{Peraire and Persson(2008)}]{CDG_scheme}
\bibinfo{author}{J.~Peraire}, \bibinfo{author}{P.-O. Persson},
\newblock \bibinfo{title}{{The Compact Discontinuous Galerkin (CDG) Method for Elliptic Problems}},
\newblock \bibinfo{journal}{SIAM Journal on Scientific Computing} \bibinfo{volume}{30} (\bibinfo{year}{2008}) \bibinfo{pages}{1806--1824}.
\bibitem[{Brezzi et~al.(2000)Brezzi, Manzini, Marini, Pietra, and Russo}]{brezzi_dg}
\bibinfo{author}{F.~Brezzi}, \bibinfo{author}{G.~Manzini}, \bibinfo{author}{D.~Marini}, \bibinfo{author}{P.~Pietra}, \bibinfo{author}{A.~Russo},
\newblock \bibinfo{title}{{Discontinuous Galerkin approximations for elliptic problems}},
\newblock \bibinfo{journal}{Numerical Methods for Partial Differential Equations: An International Journal} \bibinfo{volume}{16} (\bibinfo{year}{2000}) \bibinfo{pages}{365--378}.
\bibitem[{Shi-Dong and Nadarajah(2021)}]{doug}
\bibinfo{author}{D.~Shi-Dong}, \bibinfo{author}{S.~Nadarajah},
\newblock \bibinfo{title}{Full-space approach to aerodynamic shape optimization},
\newblock \bibinfo{journal}{Computers \& Fluids} \bibinfo{volume}{218} (\bibinfo{year}{2021}) \bibinfo{pages}{104843}.
\bibitem[{Arndt et~al.(2021)Arndt, Bangerth, Davydov, Heister, Heltai, Kronbichler, Maier, Pelteret, Turcksin, and Wells}]{dealii}
\bibinfo{author}{D.~Arndt}, \bibinfo{author}{W.~Bangerth}, \bibinfo{author}{D.~Davydov}, \bibinfo{author}{T.~Heister}, \bibinfo{author}{L.~Heltai}, \bibinfo{author}{M.~Kronbichler}, \bibinfo{author}{M.~Maier}, \bibinfo{author}{J.-P. Pelteret}, \bibinfo{author}{B.~Turcksin}, \bibinfo{author}{D.~Wells},
\newblock \bibinfo{title}{The {deal.II} finite element library: Design, features, and insights},
\newblock \bibinfo{journal}{Computers \& Mathematics with Applications} \bibinfo{volume}{81} (\bibinfo{year}{2021}) \bibinfo{pages}{407--422}.
\bibitem[{Chandrashekar(2013)}]{CH_flux}
\bibinfo{author}{P.~Chandrashekar},
\newblock \bibinfo{title}{{Kinetic Energy Preserving and Entropy Stable Finite Volume Schemes for Compressible Euler and Navier-Stokes Equations}},
\newblock \bibinfo{journal}{Communications in Computational Physics} \bibinfo{volume}{14} (\bibinfo{year}{2013}) \bibinfo{pages}{1252--1286}.
\bibitem[{Sagebaum et~al.(2019)Sagebaum, Albring, and Gauger}]{codipack}
\bibinfo{author}{M.~Sagebaum}, \bibinfo{author}{T.~Albring}, \bibinfo{author}{N.~Gauger},
\newblock \bibinfo{title}{{High-Performance Derivative Computations using CoDiPack}},
\newblock \bibinfo{journal}{ACM Transactions on Mathematical Software (TOMS)} \bibinfo{volume}{45} (\bibinfo{year}{2019}) \bibinfo{pages}{1--26}.
\bibitem[{Alexander(1977)}]{dirk}
\bibinfo{author}{R.~Alexander},
\newblock \bibinfo{title}{{Diagonally Implicit Runge–Kutta Methods for Stiff O.D.E.’s}},
\newblock \bibinfo{journal}{SIAM Journal on Numerical Analysis} \bibinfo{volume}{14} (\bibinfo{year}{1977}) \bibinfo{pages}{1006--1021}.
\bibitem[{Arnold(1982)}]{ip_arnold}
\bibinfo{author}{D.~N. Arnold},
\newblock \bibinfo{title}{{An Interior Penalty Finite Element Method with Discontinuous Elements}},
\newblock \bibinfo{journal}{SIAM Journal on Numerical Analysis} \bibinfo{volume}{19} (\bibinfo{year}{1982}) \bibinfo{pages}{742--760}.

\end{thebibliography}

\end{document}